\documentclass[11pt,a4paper]{article}
\usepackage{amsfonts}
\usepackage{amsmath,amssymb}
\usepackage{mathrsfs}
\usepackage{hyperref}
\usepackage{graphicx}
\usepackage{float}

\newtheorem{thm}{Theorem}[section]

\newtheorem{lem}[thm]{Lemma}
\newtheorem{prop}[thm]{Proposition}

\numberwithin{equation}{section}\allowdisplaybreaks

\def\leq{\leqslant}

\def\leq{\leqslant}
\def\geq{\geqslant}

\def\Real{{\mathbb{R}}}

\def\F{{\mathscr{F}}}

\begin{document}

\title{\large\bf Well-Posedness for the Two-dimensional Zakharov-Kuznetsov Equation}

\author{\normalsize \bf Minjie Shan
\footnote{Corresponding author. The  author was supported by CSC, grant 201606010025.}
\footnote{Current address: Department of Mathematics, Kyoto University, Kyoto 606-8502, Japan.}\\
\footnotesize
School of Mathematical Sciences, Peking University, Beijing 100871, P.R. China\\
\footnotesize
\texttt{ E-mail address: smj@pku.edu.cn}
} \maketitle

\thispagestyle{empty}
\begin{abstract}
Considering the initial value problem for the two-dimensional Zakharov-Kuznetsov equation
\begin{align}
u_{t} + \partial_{x}\Delta u + \partial_{x} u^2 = 0,\nonumber
\end{align}
we will show its global-in-time well-posedness in $H^{s}({\Real^2})$ when $\frac{11}{13}<s<1$ via the I-method by Colliander, Keel, Staffilani, Takaoka and Tao. Also we need to use
Strichartz estimates and local well-posedness results
 for $\frac{1}{2}<s$ . Additionally, local well-posedness for the symmetrized ZK equation in $ B^\frac{1}{2}_{2,1}(\Real^2)$ is
established in the context of atomic spaces introduced by Koch and Tataru. \\


\end{abstract}

\section{Introduction and Main Results}

This article is concerned with the Cauchy problem for the two-dimensional Zakharov-Kuznetsov equation(ZK)
\begin{align}
u_{t} + \partial_{x}\Delta u + \partial_{x} u^2 = 0,\quad
u(x,y,0)=u_0(x,y), \label{ZK}
\end{align}
where  $u=u(x,y,t)$ is a real-valued function of
$(x,y)\in\Real^2$, $t\geq0$ and $\Delta$ is the Laplace operator.

The ZK equation has been introduced by V.E. Zakharov and E.A. Kuznetsov \cite{ZK74} to describe the propagation of ionic-acoustic waves in magnetized plasma. It generalizes the
 Korteweg-de Vries equation, which is spatially one dimensional. Laedke and Spatschek \cite{EWLKHS82} first derived from the basic hydrodynamic equations the two-dimensional
Zakharov-Kuznetsov equation considered here. Lannes, Linares and Saut \cite{LLS13} justified rigorously that the ZK equation is a wave limit of the Euler-Poisson system both in two and in three dimensions.

Contrary to the Korteweg-de Vries equation or the Kadomtsev-Petviashvili equation, the ZK equation is not completely integrable, but there still exists an underlying Hamiltonian structure and the following two significant invariants,

\begin{align}
M(u)(t)= \int_{\Real^2} u^2(x,y,t) dxdy = \int_{\Real^2} u^2_0(x,y) dxdy = M(u_0) \label{M()}
\end{align}

and
\begin{align}
E(u)(t)= \int_{\Real^2} \frac{1}{2} | \nabla u(x,y,t)|^2 - \frac{1}{3} u^3(x,y,t) dxdy = E(u_0).\label{E()}
\end{align}

In the three-dimensional case, the Cauchy problem associated to Zakharov-Kuznetsov equation was shown to be locally well-posed in $H^s(\Real^3)$ for
$s> \frac{9}{8}$ by Linares and Saut \cite{FLJCS82} following the ideas of Koch and Tzvetkov  \cite{KT03}. Applying a sharp maximal function estimate in the time-weighted space,
Ribaud and Vento \cite{RV12} proved the local well-posedness in the larger data spaces $H^s(\Real^3)$ for  $s>1$ as well as in Besov space  $ B^1_{2,1}(\Real^3)$. Global well-posedness for ZK equation in $H^s(\Real^3)$ for  $s>1$ was obtained by L. Molinet and D. Pilod \cite{MP15}through taking advantage of the conservation laws, doubling time and the argument of Bourgain which he used to deal with time dependent  periodic nonlinear Schr\"odinger equations (see \cite{Bourgain 2a}) . Also, it is L. Molinet and D. Pilod who originality applied this crucial tool of atomic spaces to the well-posedness for  ZK equation.

Inspired by the ideas of Kenig, Ponce and Vega \cite{KPV93}, Faminskii \cite{F95} combined the local smoothing effects with  a maximal function estimate for the linearized equation
in order to obtain the local well-posedness for the two-dimensional ZK equation in the energy space $H^1(\Real^2)$. He also proved the global well-posedness by making use of the $L^2$ and $H^1$ conserved quantities for solutions of \eqref{ZK}. Local well-posedness  result was subsequently pushed down to $ s > \frac{3}{4}$  by Linares and Pastor \cite{LP09} via optimizing the proof of Faminskii. As far as we know, the most advanced result on local well-posedness in the two-dimensional case belongs to A. Gr\"unrock and S. Herr \cite{GH14} together with L. Molinet and D. Pilod \cite{MP15}. By applying the Fourier restriction norm method and one kind of sharp Strichartz estimates, they proved the local well-posedness in $H^s(\Real^2)$ for  $ s > \frac{1}{2}$ simultaneously. But there is no paper on the global well-posedness below $H^1(\Real^2)$ up to now.

The purpose of this paper is to derive global well-posedness when $ s < 1$ by means of I-method (see Colliander, Keel, Staffilani, Takaoka and Tao \cite{CKSTT02}) and improve  the local well-posedness theory with the help of atomic spaces. Our main results are as follows.

\begin{thm}
The initial value problem \eqref{ZK} is globally well-posed in $H^s(\Real^2)$ for $ s > \frac{11}{13}$.
\end{thm}
 \emph{Remark} $1.$ The meaning of ``globally well-posed'' is that given data  $u_0 \in H^s(\Real^2)$ and any time $T>0$, there exists a unique solution to \eqref{ZK}
\begin{center}
$ u(x,y,t)\in \mathbb{C}([0,T]; H^s(\Real^2))\cap  X^{s,\frac{1}{2}+}_T$
\end{center}
which depends continuously upon  $u_0$.
\\
\\\emph{Remark} $2.$  If one could replace the increment $N^{-\frac{1}{4}}$ in $E(I_Nu)$  on the right hand of \eqref {E incresement} with $N^{-\alpha}$ for some $ \alpha > 0 $, one could repeat the almost
 conservation law argument to prove global well-posedness of \eqref{ZK} for all  $ s > \frac{3-\alpha}{3+\alpha}$ .

\begin{thm}
The initial value problem \eqref{ZK} is locally well-posed in $ B^\frac{1}{2}_{2,1}(\Real^2)$.
\end{thm}
\emph{Remark} $3.$ The scale invariant Sobolev regularity is $ s_c= -1$. From the view of embedding $B^\frac{1}{2}_{2,1}(\Real^2)\subset H^{\frac{1}{2}}(\Real^2)\subset H^{-1}(\Real^2)$, it seems natural to consider the local well-posedness in $H^{\frac{1}{2}}(\Real^2)$ instead of $B^\frac{1}{2}_{2,1}(\Real^2)$ . Nevertheless this is not obvious. Indeed, if inequality (9) with $p=2$ and $q=\infty$ held valid in \cite{GH14} (see contribution $R_2$ on page 2066 of \cite{GH14}), the proof by Gr\"unrock and Herr \cite{GH14} could work for
$H^{\frac{1}{2}}(\Real^2)$. But it corresponds to the non-admissible endpoint of the Strichartz estimates, which can not be recovered by the use of Besov space. The new ingredient of our proof is Lemma \ref{L2U2 lema} in Section 5, which is concerned with the bilinear Strichartz estimate. What's more, we make use of the modulation decompose technique by Hadac, Herr and Koch (see  \cite{HaHeKo09}) to deal with the ``Region 4'' part in the proof of  Proposition \ref{Bilinear0}. Well-posedness in $ H^s(\Real^2)$ for $ s \leq \frac{1}{2}$ is still an open problem.
\\
\\\textbf{Organization of the paper.} Section 2 contains two subsections that provide preliminaries for global and local
wellposedness respectively. In Subsection 2.1, we first recall the Bourgain space $X^{s,b}$ and corresponding Strichartz-type estimates. Then, in Subsection 2.2, we summarize the definitions and results on $U^p$ and $V^p$ spaces. In Section 3 we show the almost conservation law. This is the important constituent for the proof of Theorem 1.1, which is eventually presented in Section 4. We prove our crucial bilinear
 estimate related to the ZK equation in Section 5. Finally, Section 6 is devoted to the proof of Theorem 1.2.

 We conclude this section with the notation given.

 Let $c < 1$, $C>3$ denote  universal constants, which can be different at different places. Given $A,B \geq 0$,  $A\lesssim B$ stands for $A\leq C \cdot B$, $A\sim B$
means that $A\lesssim B$ and $B\lesssim A$. We write $A\gg B$ to mean  $A>C \cdot B$.  We will often use the notation $c+\equiv c+\epsilon$ for some  $0<\epsilon\ll 1$. Similarly, we shall write  $c++\equiv c+2\epsilon$ and $ c-\equiv c-\epsilon$. We set $\langle a \rangle := (1+a^2)^ \frac{1}{2}$ for $a\in\Real$ and fix a smooth cut-off function
 $\chi\in C^\infty_0([-2,2])$ satisfying $\chi$ is even , nonnegative,  and $\chi=1$ on $[-1,1]$.

Throughout this paper we denote spatial variables by $x,y$ and their dual Fourier variables by $\xi,\eta$. As usual, $\tau$ is the dual variable of the time $t$. $\mathbb F(u)$ or $\hat u$ will denote space-time Fourier transform of $u$, whereas  $\F_{x,y}(u)$ or ${\widehat u}^{xy}$ will denote its Fourier transform in space.
For $s\in\Real$, $I^s_x $ and $I^s_y $ denote the one-dimensional Riesz-potential operators of order $-s$ with respect to spatial variable $x$ and $y$ . We also write $\zeta=(\xi,\eta)$, $\lambda=(\xi,\eta,\tau)$ and
$\mu =\tau -\xi^3-\eta^3 $ for brevity.

 We will use the capital letters $N,M$ to denote dyadic numbers and  write $\sum_{N\geq 1}a_N:=\sum_{n\in \mathbb{N}}a_{2^n}$ , $\sum_{N\geq M}a_N:=
\sum_{n\in \mathbb{N};2^n\geq M}a_{2^n}$ for dyadic summations.

\section{Function spaces and Strichartz estimates}

We write $\psi(x) :=\chi(x)-\chi(2x)$ and $\psi_N :=\psi(N^{-1}\cdot)$. The Littlewood-Paley multipliers are defined by
$$
 \widehat{P_0 u}=\chi(2|\zeta|)\widehat{u} \ \ \   and \ \ \
 \ \ \widehat{P_Nu}=\psi_N(|\zeta|)\widehat{u} \ \ \ for \ \  N\geq 1 \ \ .
$$
Given Lebesgue space exponents $q$, $r$ and a function $F(x,y,t)$ on $\Real^3$, we write
$$
\|F\|_{L^q_tL^r_{x,y}(\Real^3)}=(\int_{\Real}(\int_{\Real^2}|F(x,y,t)|^rdxdy)^{\frac{q}{r}}dt)^{\frac{1}{q}}.
$$
This norm will be shortened to $L^q_tL^r_{x,y}$ for readability, or to $L^r(\Real^3)$ when $q=r$.

\subsection{Bourgain spaces and estimates in $X^{s,b}$}

The Sobolev space $H^s(\Real^2)$ and the Bourgain space $X^{s,b}$ related to the linear part of \eqref{ZK} are spaces of real valued temperate distributions, defined via the norms
$$
\|u_0\|_{H^s(\Real^2)}=\|\langle \zeta \rangle^s {\widehat u_0}^{xy}(\xi,\eta)\|_{L^2(\Real^2)}
$$
\begin{align}
\|u\|_{X^{s,b}}=\|\langle \zeta \rangle^s \langle \tau-\xi^3-\xi\eta^2 \rangle^b \widehat u(\xi,\eta,\tau)\|_{L^2(\Real^3)} \label{X Norm}
\end{align}
respectively. We will need truncated versions of the Bourgain norm  \eqref{X Norm},
\begin{align}
\|u\|_{X^{s,b}_\delta}=\inf \limits_{{\tilde u=u} \ on \ [0,\delta]} \|\tilde u\|_{X^{s,b}} \label{X1 Norm}.
\end{align}

We recall some well-known estimates which play important roles to obtain almost conservation laws.
\begin{lem}\label{X Lema1}
Let $ s>\frac{1}{2}$.  We have
\begin{align}
\|e^{-t\partial_{x}\Delta}u_0\|_{X^{1, \frac{1}{2}+}_\delta} \lesssim \|u_0\|_{H^1(\Real^2)},\label{estimate1}
\end{align}
\begin{align}
\|\int^t_0 e^{-(t-t')\partial_{x}\Delta}F(t')dt'\|_{X^{1, \frac{1}{2}+}_\delta} \lesssim \|F\|_{X^{1, -\frac{1}{2}+}_\delta},\label{estimate2}
\end{align}
\begin{align}
\|F\|_{X^{1, -\frac{1}{2}+}_\delta} \lesssim \delta^{0+} \|F\|_{X^{1, -\frac{1}{2}++}_\delta},\label{estimate3}
\end{align}
\begin{align}
\|\partial_{x}(uv)\|_{X^{s, -\frac{1}{2}++}} \lesssim  \|u\|_{X^{s, \frac{1}{2}+}}\|v\|_{X^{s, \frac{1}{2}+}}.  \label{estimate4}
\end{align}
\end{lem}
{\bf Proof.} Concerning \eqref{estimate1}, \eqref{estimate2} and \eqref{estimate3}, see for example Subsection 2.3 in \cite{MP15}. The last bilinear estimate \eqref{estimate4} can be found in \cite{GH14} and \cite{MP15}.
\begin{lem} \label{X Lema2}
Let $0 \leq \epsilon < \frac{1}{2}$ and  $0 \leq \theta \leq 1 $. Assume that $2 \leq p\leq 5$, $2 \leq \tilde q$, $(q,r)$ satisfy $q=\frac{6}{\theta(2+\epsilon)}$ and $r=\frac{2}{1-\theta}$. Then, we have
\begin{align}
\|u\|_{L^p(\Real^3)} \lesssim \|u\|_{X^{0, \frac{1}{2}+}}, \label{estimate15}
\end{align}
\begin{align}
\|u\|_{L^{\tilde q}_t L^2_{x,y}(\Real^3)} \lesssim \|u\|_{X^{0, \frac{1}{2}+}}, \label{estimate5}
\end{align}
\begin{align}
\|D^{\frac{\theta \epsilon}{2}}_xu\|_{L^q_t L^r_{x,y}(\Real^3)} \lesssim \|u\|_{X^{0, \frac{1}{2}+}}.\label{estimate6}
\end{align}
In particular, it follows from \eqref{estimate6} by choosing $\epsilon = \frac{1}{2}-$ and  $ \theta = 1 $ that
\begin{align}
\|D^{\frac{1}{4}-}_xu\|_{L^{\frac{12}{5}+}_t L^{\infty}_{x,y}(\Real^3)} \lesssim \|u\|_{X^{0, \frac{1}{2}+}}\label{estimate7}
\end{align}
\end{lem}
{\bf Proof.} Strichartz estimates for the unitary group $\{e^{-t\partial_{x}\Delta}\}$ imply
\begin{align}
\|e^{-t\partial_{x}\Delta}u_0\|_{L^q_t L^r_{x,y}(\Real^3)} \lesssim \|u_0\|_{L^2(\Real^2)} \ \ \ \ \ if \ \frac{3}{q}+\frac{2}{r}=1, \ q>3. \label{estimate8}
\end{align}
From Lemma 5.3 in \cite{WaHaHuGu11}, we can obtain
\begin{align}
\|u\|_{L^q_t L^r_{x,y}(\Real^3)} \lesssim \|u\|_{X^{0, \frac{1}{2}+}} \ \ \ \ \ if \ \frac{3}{q}+\frac{2}{r}=1, \ q>3. \label{estimate9}
\end{align}
Then,the estimate \eqref{estimate15} follows from interpolation between \eqref{estimate9} for $q=r=5$ and the trivial bound $\|u\|_{L^2(\Real^3)} \lesssim \|u\|_{X^{0, \frac{1}{2}+}}$; while the estimate \eqref{estimate5} follows from interpolation between \eqref{estimate9} for $q=\infty, r=2$ and  $\|u\|_{L^2(\Real^3)} \lesssim \|u\|_{X^{0, \frac{1}{2}+}}$.
The claim \eqref{estimate6} is an immediate consequence of Proposition 2.4 in \cite{LP09} which shows the same estimate for free solution.

We will need a bilinear Strichartz estimate that can be proved in the same way as Lemma 6.5 in \cite{MP15} for two dimensional case.
\begin{lem} \label{X Lema3}
Let $N_1\ll N_2$ and $u, v \in X^{0, \frac{1}{2}+}_\delta$. Then, we have
\begin{align}
 \|P_{N_1}uP_{N_2}v\|_{L^2([0,\delta]\times\Real^2)}\lesssim \frac{N_1^{\frac{1}{2}}}{N_2}\|P_{N_1}u\|_{X^{0, \frac{1}{2}+}_\delta}\|P_{N_2}v\|_{X^{0, \frac{1}{2}+}_\delta}.\label{estimate10}
\end{align}
\end{lem}

\subsection{$U^p$, $V^p$ spaces and estimates}

In this subsection we introduce some properties of $U^p$ and $V^p$  spaces \cite{HaHeKo09,KoTa05,KoTa07,KoTa12}. After Koch and Tataru first applied $U^p$ and $V^p$ spaces to
discuss dispersive estimates for principally normal pseudodifferential operators, this sort of function space which can be treated as the development of Bourgain space has
 attracted more and more attention in the field of low regularity well-posedness theory.

Let $1\leq p <\infty$ and $\mathcal{Z}$ be the set of finite partitions $-\infty= t_0 <t_1<...< t_{K-1} < t_K =\infty$.

For any $\{t_k\}^K_{k=0} \subset \mathcal{Z}$ and $\{\phi_k\}^{K-1}_{k=0} \in L^2$ with $\sum^{K-1}_{k=0} \|\phi_k\|^p_2=1$, $\phi_0=0$. We call the function $a: \mathbb{R}\to L^2$ given by
$$
a= \sum^{K}_{k=1} 1_{[t_{k-1}, t_k)} \phi_{k-1}
$$
 a $U^p$-atom.  The atomic space is
$$
U^p:=\left\{u= \sum^\infty_{j=1} \lambda_j a_j : a_j \ U^p-atom, \lambda_j \in \mathbb{C},  \sum^\infty_{j=1} |\lambda_j|<\infty  \right\}
$$
endowed with the norm
$$
\|u\|_{U^p}:= \inf \left\{\sum^\infty_{j=1} |\lambda_j| : u= \sum^\infty_{j=1} \lambda_j a_j ,a_j \ U^p-atom, \lambda_j \in \mathbb{C} \right\}.
$$

$V^p$ is the normed space of all functions $v: \mathbb{R} \to L^2$ such that $\lim_{t\to \pm \infty} v(t)$ exist and for which the norm
$$
\|v\|_{V^p} := \sup_{\{t_k\}^K_{k=0} \in \mathcal{Z}} \left( \sum^K_{k=1} \|v(t_k)-v(t_{k-1})\|^p_{L^2}\right)^{\frac{1}{p}}  \label{v norm}
$$
is finite, where we use the convention that $v(-\infty) = \lim_{t\to -\infty} v(t)$ and $v(\infty)=0$. We denote $v \in V^p_-$ when $v(-\cdot) \in V^p $. Moreover, we define the closed subspace $V^p_{rc}$ $(V^p_{-,rc})$ of all right continuous functions in $V^p$ $(V^p_-)$.

We write $S:=-\partial^3_x-\partial^3_y$ and define the associated unitary operator $e^{tS}:L^2  \to L^2$ to be the Fourier multiplier
$$
\widehat {e^{tS}u_0}^{xy}(\xi,\eta)=e^{it(\xi^3+\eta^3)}\widehat {u_0}^{xy}(\xi,\eta).
$$

We define
\begin{itemize}
\item[\rm (i)] $U^p_S := e^{\cdot S} U^p$ with norm  $ \|u\|_{U^p_S} = \|e^{-\cdot S} u \|_{U^p}$ ,
\item[\rm (ii)] $V^p_S := e^{\cdot S} V^p$ with norm  $ \|v\|_{V^p_S} = \|e^{-\cdot S} v \|_{V^p}$.
\end{itemize}

Let us define the smooth projections
$$\widehat{Q_M u}(\xi,\eta,\tau):=\psi_M(\tau)\widehat{u}(\xi,\eta,\tau),
$$
$$
 \widehat{Q^S_M u}(\xi,\eta,\tau):=\psi_M(\tau-\xi^3-\eta^3)\widehat{u}(\xi,\eta,\tau),
$$
as well as $Q^S_{\geq M}:=\sum_{N\geq M}Q^S_N$ and $Q^S_{< M}:=I-Q^S_{\geq M}$. Note that we have
$$Q^S_M=e^{\cdot S}Q_Me^{-\cdot S}$$
and  similarly for $Q^S_{\geq M}$ and $Q^S_{< M}$.\vspace{3mm}

Here are some results in $U^p$ and $V^p$.
\begin{prop}  \label{UV prop1} Let $1< p <q <\infty$ and $\frac{1}{p}+\frac{1}{p'}=1$. We have
\begin{itemize}
 \item[\rm (i)] $U^p$, $V^p$, $V^p_{rc}$, $V^p_{-}$ and $V^p_{-,rc}$ are Banach spaces,
\item[\rm (ii)] $U^p\subset V^p_{-,rc} \subset U^q$ ,
\item[\rm (iii)]$\|u\|_{U^p}=\sup \limits_{\|v\|_{V^{p'}}=1} |\int \langle u'(t), v(t)\rangle dt|$  if $u\in V^1_{-}(\subset U^p)$ is absolutely continuous on compact interval.
\end{itemize}
\end{prop}
\begin{lem}\label{UV lema1}
 We have
\begin{align}
  \|Q^S_{\geq M}u \|_{L^2(\Real^3) }  \lesssim M^{-\frac{1}{2}} \|u\|_{V^2_S}. \label{UV estimate 1}
\end{align}
\begin{align}
  \|Q^S_{\geq M}u \|_{U^p_S} \lesssim \|u\|_{U^p_S} , \ \ \  \|Q^S_{< M}u \|_{U^p_S} \lesssim \|u\|_{U^p_S}. \label{UV estimate 2}
\end{align}
\begin{align}
 \|Q^S_{\geq M}u \|_{V^p_S} \lesssim \|u\|_{V^p_S} , \ \ \  \|Q^S_{< M}u \|_{V^p_S} \lesssim \|u\|_{V^p_S}.\label{UV estimate 3}
\end{align}
\end{lem}

Similarly to  Lemma 2.3 in \cite{GTV97} and Lemma 5.3 in  \cite{WaHaHuGu11}, we have the following general extension principle for $U^p_S$ spaces.
\begin{prop}  \label{UV prop2}
 Let $ T_0: L^2 \times \cdots \times L^2 \to L^1_{loc}(\Real^2; \mathbb{C}) $ be a n-linear operator. Assume that for some $1< p,q <\infty$
$$\| T_0(e^{\cdot S}\phi_1, \cdots, e^{\cdot S}\phi_n )\|_{L^p_t(\Real; L^q_{x,y}(\Real^2))} \lesssim \prod_{j=1}^n\|\phi_j\|_{L^2}.$$
Then, there exists $ T: U^p_S \times \cdots \times U^p_S \to L^p_t(\Real; L^q_{x,y}(\Real^2)) $ satisfying
$$\| T(u_1, \cdots, u_n )\|_{L^p_t(\Real; L^q_{x,y}(\Real^2))} \lesssim \prod_{j=1}^n\|u_j\|_{U^p_S},$$
such that $ T(u_1, \cdots, u_n )(t)(x,y)=T_0(u_1(t), \cdots, u_n(t) )(x,y) \ a.e.$.
\end{prop}

\section{Almost conservation law}

Even though $\|u\|_{H^s}$ is no longer an invariant, one modified version of the solution \eqref{ZK} has a finite energy which is almost conserved in time. In other words we need to
find some quantity that has similar properties as $E(u)(t)$. What's more, this quantity increase much slower than $\|u\|_{H^s}$ during the evolution.

Given $m:\Real^{2k} \to \mathbb{C}$, m is said to be symmetric if
$$m(\zeta_1, \cdots,\zeta_k)=m(\sigma(\zeta_1, \cdots,\zeta_k))$$
for all $\sigma \in S_k$, where $S_k$ is the permutation group for $k$ elements. We define the symmetrization of $m$ as following
$$[m]_{sym}(\zeta_1, \cdots,\zeta_k)=\frac{1}{k!}\sum_{\sigma \in S_k}m(\sigma(\zeta_1, \cdots,\zeta_k)).$$
For each $m$, a $k$-linear functional acting on $k$ functions $u_1, \cdots,u_k$ is given by
$$\Lambda_k(m;u_1, \cdots,u_k)=\int_{\zeta_1+\cdots+\zeta_k=0}m(\zeta_1, \cdots,\zeta_k)\widehat u^{xy}_1(\zeta_k)\cdots\widehat u^{xy}_k(\zeta_k).$$
We usually write $\Lambda_k(m)$ instead of $\Lambda_k(m;u, \cdots,u)$ for convenience. It is easy to see that we have  $\Lambda_k(m)= \Lambda_k([m]_{sym})$ by symmetries.

Next we introduce I-method.

Given $s<1$, $N\gg1$ and a smooth, radially symmetric, non-increasing function $m^s_N(\zeta)$ satisfying
$$ m^s_N(\zeta)=
\left\{
\begin{array}{l}
1 \ \ \ \ \ \ \ \ \ \ \ \ |\zeta|\leq N \\
(\frac{|\zeta|}{N})^{s-1} \ \ \ \ |\zeta|\geq 2N
\end{array}
\right. ,
$$
we define the Fourier multiplier operator
$$\widehat{I^s_Nf}^{xy}(\zeta):=m^s_N(\zeta) \widehat f^{xy}(\zeta).
$$
We will drop the $N$ and $s$ from the notation sometimes, writing $I$ and $m$ for simplicity.

The following general interpolation result is useful in low regularity global well-posedness theory(cf. Lemma 12.1 in \cite{CKSTT04}).
\begin{lem}\label{L5}
Let $n\geq1$. Suppose that $Z$,$X_1,\cdots, X_n$ are translation invariant Banach spaces and $T$ is a translation invariant $n$-linear operator such that
$$\|I^s_1T(u_1, \cdots,u_n)\|_Z\lesssim\prod_{j=1}^n\|I^s_1u_j\|_{X_j}
$$
for all $u_1, \cdots,u_n$ and all $\frac{1}{2}\leq s \leq 1$. Then we have
$$\|I^s_NT(u_1, \cdots,u_n)\|_Z\lesssim\prod_{j=1}^n\|I^s_Nu_j\|_{X_j}
$$
for  all $u_1, \cdots,u_n$ , all $\frac{1}{2}\leq s \leq 1$ and $N\geq 1$, with the implicit constant independent of $N$.
\end{lem}

Note that $X^{s,b}$ are translation invariant Banach spaces and $\partial_x$ is a translation invariant multi-linear operator, this interpolation lemma is available later.

Now we want to control  $\|u\|_{H^s}$ by $E(Iu)(t)$(see  \eqref {E()}). Actually these two quantities are comparable.
\begin{prop} Let $\frac{11}{13}<s<1$. We have
\begin{align}
|E(Iu)(t)|\lesssim N^{2(1-s)}\|u(t)\|^2_{\dot H^s(\Real^2)}+\|u(t)\|^3_{L^3(\Real^2)} \label{E(Iu)control}
\end{align}
\begin{align}
\|u(t)\|^2_{ H^s(\Real^2)}\lesssim |E(Iu)(t)|+\|u_0\|^2_{L^2(\Real^2)}+\|u_0\|^4_{L^2(\Real^2)}. \label{u}
\end{align}
\end{prop}
{\bf Proof.} From the fact that $|\zeta|m^s_N(\zeta)\lesssim |\zeta|^sN^{1-s} $, one has
\begin{align}
\|\nabla Iu\|^2_{L^2(\Real^2)} &  =\| |\zeta| m^s_N(\zeta) {\widehat u}^{xy}\|^2_{L^2(\Real^2)}    \notag \\
&  \lesssim\ N^{2(1-s)}\|u(t)\|^2_{\dot H^s(\Real^2)}.   \nonumber
\end{align}

$\|Iu(t)\|^3_{L^3(\Real^2)}\lesssim\|u(t)\|^3_{L^3(\Real^2)}$ follows from  H\"ormander-Mikhlin multiplier theorem, hence the energy is bounded by the right hand side of
\eqref{E(Iu)control}.

What's more, by using  the definition of $Iu$ and $L^2$ conservation \eqref {M()} , we have
\begin{align}
\|u(t)\|^2_{ H^s(\Real^2)} & \lesssim\|u(t)\|^2_{\dot H^s(\Real^2)}+\|u\|^2_{L^2(\Real^2)}     \notag \\
&  \lesssim\|\nabla Iu(t)\|^2_{L^2(\Real^2)}+\|u_0\|^2_{L^2(\Real^2)}.   \label {E3-1}
\end{align}
According to Gagliardo-Nirenberg's inequality, it holds
\begin{align}
\|\nabla Iu(t)\|^2_{L^2(\Real^2)} &\lesssim|E(Iu)(t)|+\|Iu(t)\|^3_{L^3(\Real^2)}  \notag \\
& \lesssim |E(Iu)(t)|+\|\nabla Iu(t)\|_{L^2(\Real^2)} \|Iu(t)\|^2_{L^2(\Real^2)}     \notag \\
&  \lesssim |E(Iu)(t)|+\epsilon\|\nabla Iu(t)\|^2_{L^2(\Real^2)} +C(\epsilon)\|Iu(t)\|^4_{L^2(\Real^2)} \notag \\
&  \lesssim|E(Iu)(t)|+\epsilon\|\nabla Iu(t)\|^2_{L^2(\Real^2)} +C(\epsilon)\|u_0\|^4_{L^2(\Real^2)} . \nonumber
\end{align}
Let $\epsilon\ll 1$, one can obtain
\begin{align}
\|\nabla Iu(t)\|^2_{L^2(\Real^2)} &\lesssim|E(Iu)(t)|+\|u_0\|^4_{L^2(\Real^2)}.  \label {E3-2}
\end{align}
Then, \eqref {u} follows from \eqref {E3-1} and \eqref {E3-2}.

We need to control for small times the smoothed solution before extending to a global one. Here is a modified local existence theorem.
\begin{prop} \label{p2}
Let $\frac{11}{13}<s<1$. Assume $u_0$ satisfies $|E(Iu_0)|\leq 1$. Then there is a constant $\delta=\delta(\|u_0\|_{L^2(\Real^2)})$ and a unique solution $u$ to \eqref {ZK} on $[0,\delta]$,
 such that
$$\|Iu\|_{X^{1, \frac{1}{2}+}_\delta}\lesssim 1
$$
where the implicit constant is independent of $\delta$.
\end{prop}
{\bf Proof.} We mimic the usual iteration argument showing local well-posedness.

Acting multiplier operator $I$ on both sides of \eqref {ZK}, one can obtain
\begin{align}\partial_{t}Iu + \partial_{x}\Delta Iu + \partial_{x} I(u^2) = 0.\label{equationI}\end{align}
We rewrite the differential equation as an integral equation by Duhamel's principle
$$Iu=e^{-t\partial_{x}\Delta}Iu_0-\int^t_0e^{-(t-t')\partial_{x}\Delta}\partial_{x}I(u^2)(t')dt'.$$
Estimates \eqref {estimate1}-\eqref {estimate3} give us
\begin{align}
\|Iu\|_{X^{1, \frac{1}{2}+}_\delta} &\lesssim\|e^{-t\partial_{x}\Delta}Iu_0\|_{X^{1, \frac{1}{2}+}_\delta}+\|\int^t_0e^{-(t-t')\partial_{x}\Delta}\partial_{x}I(u^2)(t')dt'\|_{X^{1, \frac{1}{2}+}_\delta}  \notag \\
& \lesssim \|Iu_0\|_{H^1(\Real^2)} +\|\partial_{x}I(u^2)\|_{X^{1, -\frac{1}{2}+}_\delta}    \notag \\
&   \lesssim \|Iu_0\|_{H^1(\Real^2)} +\delta^{0+}\|\partial_{x}I(u^2)\|_{X^{1, -\frac{1}{2}++}_\delta}. \label{3.5}
\end{align}
By the definition of the restricted norm \eqref{X1 Norm}, we can choose $\tilde u \in X^{1, \frac{1}{2}+} $ satisfies $\tilde u|_{[0, \delta]}=u$,
\begin{align}\|Iu\|_{X^{1, \frac{1}{2}+}_\delta} \sim \|I\tilde u\|_{X^{1, \frac{1}{2}+}}\end{align} \label{3.6}
and
\begin{align}
\|\partial_{x}I(u^2)\|_{X^{1, -\frac{1}{2}++}_\delta}\lesssim \|\partial_{x}I(\tilde u^2)\|_{X^{1, -\frac{1}{2}++}}. \label{3.7}
\end{align}
We will show shortly that
\begin{align}
\|\partial_{x}I(\tilde u^2)\|_{X^{1, -\frac{1}{2}++}}\lesssim  \|I\tilde u\|^2_{X^{1, \frac{1}{2}+}}.\label{3.8}
\end{align}
Using the Lemma \ref{L5}, it suffices to prove
\begin{align}
\|I^s_1\partial_{x}(\tilde u^2)\|_{X^{1, -\frac{1}{2}++}}\lesssim  \|I^s_1\tilde u\|^2_{X^{1, \frac{1}{2}+}}. \label{3.9}
\end{align}
Note that $\|I^s_1F\|_{X^{1, b}} \sim \|F\|_{X^{s, b}}$, \eqref{3.9} is an immediate consequence of \eqref{estimate4}.

Combining \eqref{3.5}-\eqref{3.8}, we have
\begin{align}\|Iu\|_{X^{1, \frac{1}{2}+}_\delta}\lesssim \|Iu_0\|_{H^1(\Real^2)} +\delta^{0+}\|Iu\|^2_{X^{1, \frac{1}{2}+}_\delta},  \label{3.10} \end{align}
and similarly
$$\|Iu-Iv\|_{X^{1, \frac{1}{2}+}_\delta}\lesssim \delta^{0+}(\|Iu\|_{X^{1, \frac{1}{2}+}_\delta}+\|Iv\|_{X^{1, \frac{1}{2}+}_\delta})\|Iu-Iv\|_{X^{1, \frac{1}{2}+}_\delta}.$$
Then, one can obtain the local well-posedness by means of the contraction mapping principle.

Moreover, setting $Q(\delta)\equiv \|Iu\|_{X^{1, \frac{1}{2}+}_\delta}$, the bound \eqref{3.10} yields
$$Q(\delta)\lesssim \|Iu_0\|_{H^1(\Real^2)} +\delta^{0+}(Q(\delta))^2. $$
The proof of \eqref {u} gives us
\begin{align} \|Iu_0\|_{H^1(\Real^2)} & \lesssim|E(Iu_0)|^{\frac{1}{2}} +\|u_0\|_{L^2(\Real^2)} +\|u_0\|^2_{L^2(\Real^2)}  \notag \\
&   \lesssim 1+\|u_0\|_{L^2(\Real^2)} +\|u_0\|^2_{L^2(\Real^2)}\nonumber .
\end{align}
As $Q(\delta)$ is continuous in the variable  $\delta$, a bootstrap argument yields $\|Iu\|_{X^{1, \frac{1}{2}+}_\delta}\lesssim 1.$

We consider the growth of $E(Iu)(t)$.

Using the definition of $E(Iu)$ , equation \eqref {equationI} and integration by parts,
\begin{align}
\frac{dE(Iu)(t)}{dt}& =\int_{\Real^2}[\nabla Iu \cdot \nabla I\partial_{t}u-(Iu)^2I\partial_{t}u ]dxdy \notag \\
& =- \int_{\Real^2}(I\partial_{t}u) [\Delta Iu +(Iu)^2 ]dxdy   \notag \\
& = \int_{\Real^2}\partial_{x}(I\Delta u+Iu^2)  [\Delta Iu +(Iu)^2 ]dxdy   \notag \\
& = \int_{\Real^2}(\partial_{x}I\Delta u)[(Iu)^2-Iu^2]+(\partial_{x}Iu^2)(Iu)^2dxdy.   \nonumber
\end{align}
When integrating in time and applying the Parseval formula, one has
\begin{align}
E(Iu)(\delta)- E(Iu)(0)& =-i\int^{\delta}_0\int_{\sum^3_{j=1}\zeta_j=0}\xi_1|\zeta_1|^2[1-\frac{m(\zeta_2+\zeta_3)}{m(\zeta_2)m(\zeta_3)}]\prod^3_{j=1}\widehat{Iu}^{xy}(\zeta_j)
 \notag \\
& +i\int^{\delta}_0\int_{\sum^4_{j=1}\zeta_j=0}(\xi_1+\xi_2)\frac{m(\zeta_1+\zeta_2)}{m(\zeta_1)m(\zeta_2)}\prod^4_{j=1}\widehat{Iu}^{xy}(\zeta_j). \label{In1}
\end{align}

We estimate these two terms on the right hand of \eqref {In1} respectively.
\begin{prop} \label{p3}
Let $\frac{11}{13}<s<1$. We have
\begin{align}
 |\int^{\delta}_0 \Lambda_3( \xi_1|\zeta_1|^2[1-\frac{m(\zeta_2+\zeta_3)}{m(\zeta_2)m(\zeta_3)}];Iu)|\lesssim
  N^{-\frac{1}{4}+}\|Iu\|^3_{X^{1,\frac{1}{2}+}_{\delta}} .\label{Lambda(3)}
\end{align}
\end{prop}
{\bf Proof.}  We may assume $\widehat u^{xy}$ is non-negative. Firstly, we break $u$ into a sum of dyadic constituents $P_{N_k}u$, each with frequency support $|\zeta|\sim N_k$, $k=0,1,\cdots$. Once we show
\begin{align}
| & \int^{\delta}_0\int_{\sum^3_{j=1}\zeta_j=0}[\xi_1|\zeta_1|^2(1-\frac{m(\zeta_2+\zeta_3)}{m(\zeta_2)m(\zeta_3)})]_{sym}\prod^3_{j=1}\widehat{u_j}^{xy}(\zeta_j)| \notag \\
&\lesssim N^{-\frac{1}{4}+}(N_1N_2N_3)^{0-}\prod^3_{j=1}\|u_j\|_{X^{1,\frac{1}{2}+}_{\delta}} \label{Lambda(3.1)}
\end{align}
 for any function $u_j$, $j=1,2,3$ with frequencies supported on $|\zeta_j|\sim  N_j$, we conclude our desired bound \eqref {Lambda(3)} by summing over all dyadic pieces $P_{N_k}u$.

We denote $T_1$ the left hand of \eqref{Lambda(3.1)} and $M(\zeta_1,\zeta_2,\zeta_3)=[\xi_1|\zeta_1|^2(1-\frac{m(\zeta_2+\zeta_3)}{m(\zeta_2)m(\zeta_3)})]_{sym}$. One can  assume $N_1\sim N_2\geq N_3$ by symmetry.\\
{\bf Case 1.} $N\gg N_1$

In this case, the symbol $M(\zeta_1,\zeta_2,\zeta_3)$ is identically zero and the bound holds trivially.\\
{\bf Case 2.} $N_1\sim N_2\sim N_3\gtrsim N$

Since $m(\zeta_1)\sim m(\zeta_2)\geq m(\zeta_2)m(\zeta_3)$, then
$$|1-\frac{m(\zeta_2+\zeta_3)}{m(\zeta_2)m(\zeta_3)}|\lesssim\frac{m(\zeta_1)}{m(\zeta_2)m(\zeta_3)}\lesssim (\frac{N_1}{N})^{1-s},$$
hence
$$|M(\zeta_1,\zeta_2,\zeta_3)|\lesssim (\frac{N_1}{N})^{1-s} N_1^{\frac{11}{4}+}(|\xi_1|^{\frac{1}{4}-}+|\xi_2|^{\frac{1}{4}-}+|\xi_3|^{\frac{1}{4}-}).$$
Applying H\"{o}rder's inequalities, \eqref{estimate5} and \eqref{estimate7}, we obtain
\begin{align}
  T_1 & \lesssim(\frac{N_1}{N})^{1-s} N_1^{\frac{11}{4}+}\|\partial^{\frac{1}{4}-}_x u_1\|_{L^{\frac{12}{5}+}_{[0,\delta]}L^\infty_{x,y}}
\|u_2\|_{L^{\frac{24}{7}-}_{[0,\delta]}L^2_{x,y}}\|u_3\|_{L^{\frac{24}{7}-}_{[0,\delta]}L^2_{x,y}} \notag \\
&\lesssim(\frac{N_1}{N})^{1-s} N_1^{\frac{11}{4}+}\prod^3_{j=1}\|u_j\|_{X^{0,\frac{1}{2}+}_\delta}\notag \\
&\lesssim N^{s-1} N_1^{\frac{3}{4}-s+}\prod^3_{j=1}\|u_j\|_{X^{1,\frac{1}{2}+}_\delta} \notag \\
&\lesssim N^{-\frac{1}{4}+} N_1^{0-}\prod^3_{j=1}\|u_j\|_{X^{1,\frac{1}{2}+}_\delta}. \nonumber
\end{align}
{\bf Case 3.} $N_1\sim N_2\gg N_3, N_1 \gtrsim N $

We would like to use mean value formulas to get some cancelation for the symbol.
We write $M=M_1-M_2$, where $M_1=\sum^3_{j=1}\xi_j|\zeta_j|^2$ and
$M_2=\frac{\sum^3_{j=1}\xi_j|\zeta_j|^2m^2(\zeta_j)}{m(\zeta_1)m(\zeta_2)m(\zeta_3)}$.
The corresponding terms are
\begin{align}T_{1,1}=|\int^{\delta}_0\int_{\sum^3_{j=1}\zeta_j=0}M_1(\zeta_1,\zeta_2,\zeta_3)\prod^3_{j=1}\widehat{u_j}^{xy}(\zeta_j)| \nonumber \end{align}
and \begin{align}T_{1,2}=|\int^{\delta}_0\int_{\sum^3_{j=1}\zeta_j=0}M_2(\zeta_1,\zeta_2,\zeta_3)\prod^3_{j=1}\widehat{u_j}^{xy}(\zeta_j)|. \nonumber \end{align}
It is obvious that $T_1\leq T_{1,1}+T_{1,2}$.

First we consider the contribution of $T_{1,1}$.
On the hyperplane $\left\{(\zeta_1,\zeta_2,\zeta_3):\zeta_1+\zeta_2+\zeta_3=0\right\}$, we have
\begin{align}M_1&=\xi_1|\zeta_1|^2-(\xi_1+\xi_3)|\zeta_1+\zeta_3|^2+\xi_3|\zeta_3|^2\notag\\
&=\xi_1(|\zeta_1|^2-|\zeta_1+\zeta_3|^2)+\xi_3(|\zeta_3|^2-|\zeta_2|^2).\nonumber \end{align}
Since $|\zeta_3|\ll |\zeta_1|$, mean-value theorem for vector-valued functions helps to bound $M_1$,
\begin{align}|M_1|\lesssim |\xi_1||\zeta_1||\zeta_3|+|\xi_3||\zeta_2|^2
\lesssim N_1^2N_3.\nonumber \end{align}
Hence,
\begin{align}
T_{1,1}&\lesssim N_1^2N_3\|u_1\|_{L^2([0,\delta]\times \Real^2)}\|u_2u_3\|_{L^2([0,\delta]\times \Real^2)}
\notag \\
&\lesssim N_1^2N_3\frac{N_3^{\frac{1}{2}}}{N_1}\prod^3_{j=1}\|u_j\|_{X^{0,\frac{1}{2}+}_\delta}\notag \\
&\lesssim N^{-\frac{1}{2}+}N_1^{0-}\prod^3_{j=1}\|u_j\|_{X^{1,\frac{1}{2}+}_\delta},
\nonumber
 \end{align}
the second inequality is due to Lemma \ref{X Lema3}.

In the next place we consider the contribution of $T_{1,2}$.
With regard to the symbol $M_2$, we have
\begin{align}\sum^3_{j=1}\xi_j|\zeta_j|^2m^2(\zeta_j)&=\xi_1|\zeta_1|^2m^2(\zeta_1) -(\xi_1+\xi_3)|\zeta_1+\zeta_3|^2m^2(\zeta_1+\zeta_3)+\xi_3|\zeta_3|^2m^2(\zeta_3)\notag \\
&=\xi_1[|\zeta_1|^2m^2(\zeta_1)-|\zeta_1+\zeta_3|^2m^2(\zeta_1+\zeta_3)]\notag \\
&+\xi_3[|\zeta_3|^2m^2(\zeta_3)-|\zeta_2|^2m^2(\zeta_2)].
\nonumber \end{align}
Mean-value theorem tells us that
$$||\zeta_1|^2m^2(\zeta_1)-|\zeta_1+\zeta_3|^2m^2(\zeta_1+\zeta_3)|\lesssim |\zeta_1||\zeta_3|m^2(\zeta_1).$$
For the other part, we need split the different frequency interactions into two subcases according to the size of the parameter $N$ in comparision to $N_3$.
\\{\bf Case 3(a).}  $N_1\sim N_2\gg N_3\gtrsim N$

In this case, one has
$$|\zeta_3|^2m^2(\zeta_3)\ll |\zeta_1|^2m^2(\zeta_1),$$
thus we can control the symbol by
$$|M_2|\lesssim\frac{N_1^2N_3m^2(\zeta_1)}{m(\zeta_3)m^2(\zeta_1)}\lesssim N_1^2N_3(\frac{N_3}{N})^{1-s}.$$
Then, by \eqref {estimate10} we have
\begin{align}T_{1,2}&\lesssim N_1^2N_3(\frac{N_3}{N})^{1-s}\|u_1\|_{L^2([0,\delta]\times \Real^2)}\|u_2u_3\|_{L^2([0,\delta]\times \Real^2)}\notag \\
&\lesssim N_1^2N_3(\frac{N_3}{N})^{1-s}\frac{N_3^{\frac{1}{2}}}{N_1}\prod^3_{j=1}\|u_j\|_{X^{0,\frac{1}{2}+}_\delta}\notag \\
&\lesssim N^{s-1}N_1^{-1}N_3^{\frac{3}{2}-s}\prod^3_{j=1}\|u_j\|_{X^{1,\frac{1}{2}+}_\delta}\notag \\
&\lesssim N^{s-1}N_1^{0-}N_3^{\frac{1}{2}-s+}\prod^3_{j=1}\|u_j\|_{X^{1,\frac{1}{2}+}_\delta}\notag \\
&\lesssim N^{-\frac{1}{2}+}N_1^{0-}\prod^3_{j=1}\|u_j\|_{X^{1,\frac{1}{2}+}_\delta}.\nonumber \end{align}
{\bf Case 3(b).}  $N_1\sim N_2\gtrsim N\gg N_3$

Now $m(\zeta_3)=1$, therefore one has
\begin{align}|M_2|&\lesssim\frac{N_1^2N_3m^2(\zeta_1)+N_3^3}{m^2(\zeta_1)}\notag \\
&\lesssim N_1^2N_3+N_3^3(\frac{N_1}{N})^{2(1-s)}.\nonumber \end{align}
We obtain the bound of $T_{1,2}$,
\begin{align}T_{1,2}&\lesssim [N_1^2N_3+N_3^3(\frac{N_1}{N})^{2(1-s)}]\|u_1\|_{L^2([0,\delta]\times \Real^2)}\|u_2u_3\|_{L^2([0,\delta]\times \Real^2)}\notag \\
&\lesssim [N_1^2N_3+N_3^3(\frac{N_1}{N})^{2(1-s)}]\frac{N_3^{\frac{1}{2}}}{N_1}\prod^3_{j=1}\|u_j\|_{X^{0,\frac{1}{2}+}_\delta}\notag \\
&\lesssim (N_1^{-1}N_3^{\frac{1}{2}}+N^{2(s-1)}N_1^{-1-2s}N_3^{\frac{5}{2}})\prod^3_{j=1}\|u_j\|_{X^{1,\frac{1}{2}+}_\delta}\notag \\
&\lesssim (N_1^{-\frac{1}{2}}+N^{2s+\frac{1}{2}}N_1^{-1-2s})\prod^3_{j=1}\|u_j\|_{X^{1,\frac{1}{2}+}_\delta}\notag \\
&\lesssim (N^{-\frac{1}{2}+}+N^{2s+\frac{1}{2}}N^{-1-2s+})N_1^{0-}\prod^3_{j=1}\|u_j\|_{X^{1,\frac{1}{2}+}_\delta}\notag \\
&\lesssim N^{-\frac{1}{2}+}N_1^{0-}\prod^3_{j=1}\|u_j\|_{X^{1,\frac{1}{2}+}_\delta}.\nonumber \end{align}
This complete the proof of \eqref{Lambda(3.1)}, and hence \eqref{Lambda(3)}.
\begin{prop} \label{p4}
Let $\frac{11}{13}<s<1$. We have
\begin{align}
 |\int^{\delta}_0 \Lambda_4((\xi_1+\xi_2)\frac{m(\zeta_1+\zeta_2)}{m(\zeta_1)m(\zeta_2)};Iu)|\lesssim N^{-1+}\|Iu\|^4_{X^{1,\frac{1}{2}+}_{\delta}}\label{Lambda(4)}.
\end{align}
\end{prop}
{\bf Proof.} As our previous discussion of $\Lambda_3$, it suffices to show
\begin{align}
| & \int^{\delta}_0\int_{\sum^4_{j=1}\zeta_j=0}(\xi_1+\xi_2)\frac{m(\zeta_1+\zeta_2)}{m(\zeta_1)m(\zeta_2)}\prod^4_{j=1}\widehat u_j^{xy}(\zeta_j)| \notag \\
&\lesssim N^{-1+}\prod^4_{j=1}N_j^{-\alpha_j}\|u_j\|_{X^{1,\frac{1}{2}+}_{\delta}} \label{Lambda(4.1)}
\end{align}
 for any function $u_j$, $j=1,2,3,4$ with non-negative spatial frequencies supported on $|\zeta_j|\sim  N_j$ and some $\alpha_j > 0$, $j=1,2,3,4$.

We can assume at least one of $N_1,N_2,N_3,N_4$ is not smaller than $N$. Otherwise, it is easy to know
 $$[(\xi_1+\xi_2)\frac{m(\zeta_1+\zeta_2)}{m(\zeta_1)m(\zeta_2)}]_{sym}=0,$$
 the bound \eqref {Lambda(4.1)} holds trivially.
 Furthermore, one may assume $N_1\geq N_2, N_3\geq N_4$ by the symmetry of the multiplier.

We  denote $T_2$ the left hand of \eqref{Lambda(4.1)} and divide the interactions into three subcases depending on relationships between $N_1$ and $N_3$.
\\{\bf Case 1.} $N_1\gtrsim N_3, N_1 \sim N_2\gtrsim N$

We use in this instance a pointwise bound on the symbol,
$$|(\xi_1+\xi_2)\frac{m(\zeta_1+\zeta_2)}{m(\zeta_1)m(\zeta_2)}|\lesssim N_1(\frac{N_1}{N})^{2(1-s)}.$$
From \eqref {estimate15}, we have
\begin{align}
T_2&\lesssim N_1^{3-2s}N^{2(s-1)}\prod^4_{j=1}\|u_j\|_{L^4([0,\delta]\times \Real^2)}\notag \\
&\lesssim N_1^{3-2s}N^{2(s-1)}\prod^4_{j=1}\|u_j\|_{X^{0,\frac{1}{2}+}_{\delta}}\notag \\
&\lesssim N_1^{3-2s}N^{2(s-1)}\prod^4_{j=1}N_j^{-1}\|u_j\|_{X^{1,\frac{1}{2}+}_{\delta}}\notag \\
&\lesssim N^{-1+}N_1^{0-}N_3^{-1}N_4^{-1}\prod^4_{j=1}\|u_j\|_{X^{1,\frac{1}{2}+}_{\delta}}. \nonumber
\end{align}
{\bf Case 2.} $N_1\gtrsim N_3, N_1 \gg N_2$

In this case, $N_1 \sim N_3\gtrsim N$, it is easy to see that
$$|(\xi_1+\xi_2)\frac{m(\zeta_1+\zeta_2)}{m(\zeta_1)m(\zeta_2)}|\lesssim N_1(\frac{N_2}{N})^{1-s}.$$
We control $T_2$ by using \eqref {estimate15} and \eqref{estimate10},
\begin{align}
T_2&\lesssim N^{s-1}N_1N_2^{1-s}\|u_1u_2\|_{L^2([0,\delta]\times \Real^2)}\|u_3\|_{L^4([0,\delta]\times \Real^2)}\|u_4\|_{L^4([0,\delta]\times \Real^2)}\notag \\
&\lesssim N^{s-1}N_1N_2^{1-s}\frac{N_2^{\frac{1}{2}}}{N_1}\prod^4_{j=1}\|u_j\|_{X^{0,\frac{1}{2}+}_{\delta}}\notag \\
&\lesssim N^{s-1}N_1^{-2}N_2^{\frac{1}{2}-s}N_4^{-1}\prod^4_{j=1}\|u_j\|_{X^{1,\frac{1}{2}+}_{\delta}}\notag \\
&\lesssim N^{-(3-s)+}N_1^{0-}N_2^{\frac{1}{2}-s}N_4^{-1}\prod^4_{j=1}\|u_j\|_{X^{1,\frac{1}{2}+}_{\delta}}. \nonumber
\end{align}
{\bf Case 3.} $N_1\ll N_3$

That's to say $N_3\sim N_4\gg N_1\geq N_2, N_3\sim N_4\gtrsim N.$ Hence, we have
$$|(\xi_1+\xi_2)\frac{m(\zeta_1+\zeta_2)}{m(\zeta_1)m(\zeta_2)}|\lesssim N_1(\frac{N_1}{N})^{1-s}(\frac{N_2}{N})^{1-s}.$$
Pairing $u_1u_3$ and $u_2u_4$ in $L^2$ and applying \eqref{estimate10} again, it gives
\begin{align}T_2&\lesssim N_1(\frac{N_1}{N})^{1-s}(\frac{N_2}{N})^{1-s}\|u_1u_3\|_{L^2([0,\delta]\times \Real^2)}\|u_2u_4\|_{L^2([0,\delta]\times \Real^2)}\notag \\
&\lesssim N_1(\frac{N_1}{N})^{1-s}(\frac{N_2}{N})^{1-s}\frac{N_1^{\frac{1}{2}}}{N_3}\frac{N_2^{\frac{1}{2}}}{N_4}\prod^4_{j=1}\|u_j\|_{X^{0,\frac{1}{2}+}_{\delta}}\notag \\
&\lesssim N^{2(s-1)}N_1^{-(s-\frac{1}{2})}N_2^{-(s-\frac{1}{2})}N_3^{-3}\prod^4_{j=1}\|u_j\|_{X^{1,\frac{1}{2}+}_{\delta}}\notag \\
&\lesssim N^{-(5-2s)+}N_1^{-(s-\frac{1}{2})}N_2^{-(s-\frac{1}{2})}N_3^{0-}\prod^4_{j=1}\|u_j\|_{X^{1,\frac{1}{2}+}_{\delta}}.
\nonumber
\end{align}
This complete the proof of \eqref{Lambda(4.1)}, and hence \eqref{Lambda(4)}.
\begin{prop} \label {Prop E incresement}
Let $\frac{11}{13}<s<1, N\gg 1$. Assume $u_0$ satisfies $|E(Iu_0)|\leq 1$. Then there is a constant $\delta=\delta(\|u_0\|_{L^2(\Real^2)})>0$ so that there exists a unique solution
$$u(x,y,t)\in C([0,\delta],H^s(\Real^2))\cap X^{s,\frac{1}{2}+}_{\delta}$$
 of \eqref {ZK} satisfying
\begin{align}
E(I_Nu)(\delta)=E(I_Nu)(0)+O(N^{-\frac{1}{4}+}). \label {E incresement}
 \end{align}

\end{prop}
{\bf Proof.} We know from Proposition \ref{p2} that there is a unique solution $u$ to \eqref {ZK} on $[0,\delta]$ satisfying $\|Iu\|_{X^{1, \frac{1}{2}+}_\delta}\lesssim 1$.

It follows from \eqref {In1}, Proposition \ref{p3} and Proposition \ref {p4} that
\begin{align}
 |E(Iu)(\delta)- E(Iu)(0)|& =|\int^{\delta}_0 \Lambda_3( \xi_1|\zeta_1|^2[1-\frac{m(\zeta_2+\zeta_3)}{m(\zeta_2)m(\zeta_3)}];Iu)\notag \\
 &+\int^{\delta}_0 \Lambda_4((\xi_1+\xi_2)\frac{m(\zeta_1+\zeta_2)}{m(\zeta_1)m(\zeta_2)};Iu)|\notag \\
 & \lesssim N^{-\frac{1}{4}+}\|Iu\|^3_{X^{1,\frac{1}{2}+}_{\delta}}+N^{-1+}\|Iu\|^4_{X^{1,\frac{1}{2}+}_{\delta}}\notag \\
 & \lesssim N^{-\frac{1}{4}+}.\nonumber
\end{align}
Thus we prove this proposition.
\section{Global well-posedness}
For any given $u_0\in H^s(\Real^2)$, $\frac{11}{13}<s<1$ and time $T>0$, our purpose is to construct a solution on $[0,T]$. Note that the initial value problem \eqref {ZK} has a
scaling symmetry. That is, if $u$ is a solution to \eqref {ZK}, so is $u_\lambda(x,y,t)=\lambda^2u(\lambda x,\lambda y,\lambda^3t)$. It is easy to see that when $u_\lambda$ exists on $[0,\frac{T}{\lambda^3}]$, $u$ exists on on $[0,T]$.

Using \eqref {E(Iu)control}, the following energy can be made arbitrarily small by taking $\lambda$ small,
\begin{align}
E(I_Nu_{\lambda,0})&\leq N^{2(1-s)}\|u_{\lambda,0}\|^2_{\dot H^s(\Real^2)}+\|u_{\lambda,0}\|^3_{L^3(\Real^2)}\notag \\
 & \leq N^{2(1-s)}\lambda^{2(s+1)}\|u_0\|^2_{\dot H^s(\Real^2)}+\lambda^4\|u_0\|^3_{L^3(\Real^2)}\notag \\
 & \leq C_0(N^{2(1-s)}\lambda^{2(s+1)}+\lambda^4)(1+\|u_0\|_{H^s(\Real^2)})^3.\nonumber
\end{align}
Assuming $N\gg1$ is given ($N$ will be chose later), we choose our scaling parameter
\begin{align}
\lambda=\lambda(N,\|u_0\|_{H^s(\Real^2)})\sim N^{\frac{s-1}{s+1}}\nonumber
\end{align}
such that $E(I_Nu_{\lambda,0})\leq \frac{1}{4}$. This is feasible because $\lambda^4\sim N^{\frac{4(s-1)}{s+1}}\ll1$.

 We can now apply Proposition \ref {Prop E incresement} to the scaled initial data $u_{\lambda,0}$, then we get
$$E(Iu_\lambda)(\delta)\leq\frac{1}{4}+CN^{-\frac{1}{4}+}<\frac{1}{2}.$$
Thus the solution $u_\lambda$ can be extended to $t\in[0,2\delta]$ by Proposition \ref {p2}.

Iterating this procedure $M$ steps, we get for $t\in[0,(M+1)\delta]$
$$E(Iu_\lambda)(t)\leq\frac{1}{4}+CMN^{-\frac{1}{4}+}$$ as long as $MN^{-\frac{1}{4}+}\lesssim1$. It means that the solution $u_\lambda$ can be extended to $t\in[0,N^{\frac{1}{4}-}\delta].$
Take $N(T)$ sufficiently large so that
$$N^{\frac{1}{4}-}\delta >\frac{T}{\lambda^3}\sim TN^{\frac{3(1-s)}{s+1}}.$$
Thus,
$$T\sim N^{\frac{13s-11}{4(s+1)}-}.$$
Note that the exponent of $N$ above is positive provided $s>\frac{11}{13}$, hence the definition of $N$ makes sense for arbitrary $T$.

In the end, we give the increment property of the solution.

From scaling and \eqref{u}, one has
\begin{align}\|u(T)\|_{ H^s(\Real^2)}&\lesssim \lambda^{-s-1}\|u_\lambda(\frac{T}{\lambda^3})\|_{ H^s(\Real^2)}\notag \\
&\lesssim \lambda^{-s-1}(|E(Iu_\lambda)(\frac{T}{\lambda^3})|^{\frac{1}{2}}+\|u_{\lambda,0}\|_{L^2(\Real^2)}+\|u_{\lambda,0}\|^2_{L^2(\Real^2)})\notag \\
&\lesssim \lambda^{-s-1} (1+\|u_0\|_{ H^s(\Real^2)})^2.\nonumber\end{align}
As $\lambda\sim T^{\frac{4(s-1)}{13s-11}+}$, the global solution $u(x,y,t)$ satisfies
$$\|u(T)\|_{ H^s(\Real^2)}\lesssim T^{\frac{4(1-s)(1+s)}{13s-11}-} (1+\|u_0\|_{ H^s(\Real^2)})^2.$$

\section{Bilinear estimate}
We turn to local well-posedness for equation \eqref{ZK}. In order to symmetrize the equation, we perform a linear change of variables. Ben-Artzi, Koch and Saut studied systematically such transformations in connection with dispersive estimates for cubic phase functions of two variables in \cite{BKS03}. We also refer to Gr\"unrock and Herr's paper \cite{GH14}. After a rotation of variables, we may consider the initial value problem
\begin{align}
u_{t} + (\partial^3_{x}+ \partial^3_{y}) u + (\partial_{x}+ \partial_{y}) u^2 = 0,\quad
u(x,y,0)=u_0(x,y)\in H^s(\Real^2), \label{ZK1}
\end{align}
instead of \eqref{ZK} without changing the well-posedness theory.

Let us list some useful estimates in $U^p_S$ and $V^p_S$.
\begin{lem} \label{in lema}
Given $N_1\geq N_2$. Assume that $(q,r)$ satisfy $\frac{3}{q}+\frac{2}{r}=1$ and $q>3$. Let $I^s_{x,-}$ be the bilinear operator with symbol $|\xi_1-\xi_2|^s$, i.e.
$$\F_{x,y}(I^s_{x,-}(f_1,f_2))(\xi,\eta)=\int_{\zeta=\zeta_1+\zeta_2}|\xi_1-\xi_2|^s \prod^2_{j=1}\widehat{f_j}^{xy}(\xi_j,\eta_j).$$
Then, we have
\begin{align}
 \|\chi(\frac{t}{T})u \|_{L^2(\Real^3)}\lesssim T^{\frac{1}{2}}\|u\|_{V^2_S}, \label{estimate51}
\end{align}
\begin{align}
  \|u\|_{L^q_tL^r_{xy}}  \lesssim \|u\|_{U^q_S},\label{estimate52}
\end{align}
\begin{align}
  \|I^{\frac{1}{8}}_xI^{\frac{1}{8}}_ye^{tS}u_0\|_{L^4(\Real^3)}  \lesssim \|u_0\|_{L^2(\Real^2)},\label{estimate53}
\end{align}
\begin{align}
  \|I^{\frac{1}{2}}_xI^{\frac{1}{2}}_{x,-}(P_{N_1}e^{tS}u_0,P_{N_2}e^{tS}v_0) \|_{L^2(\Real^3)} \lesssim N_2^{\frac{1}{2}}\|P_{N_1}u_0\|_{L^2(\Real^2)}\|P_{N_2}v_0\|_{L^2(\Real^2)} , \label{estimate54}
\end{align}
where $U^p_S$ and $V^p_S$ are defined as in Subsection 2.2. In addition, \eqref{estimate54} is equally valid with $x$ replaced by $y$.
\end{lem}
{\bf Proof.} Since $\|\chi(\frac{t}{T})u \|_{L^2(\Real^3)}  \lesssim T^{\frac{1}{2}}\|u\|_{L^\infty_tL^2_{x,y}}\lesssim T^{\frac{1}{2}}\|u\|_{V^2_S}$, we prove \eqref{estimate51}.
Estimate (8) of \cite{GH14} shows that \eqref{estimate52} holds true for free solution. Thus, the claim \eqref{estimate52} follows from Proposition \ref{UV prop2}. \eqref{estimate53} and \eqref{estimate54} can also be found in \cite{GH14}.

Similar to the technique in \cite{MV15}, one can decompose the time cut-off into low- and high-frequency parts.

For any $T>0$, we write $1_T$ the characteristic function of $[0,T]$ and
$$1_T=1^{low}_{T,L}+1^{high}_{T,L}, \ \ \widehat {1^{low}_{T,L}}(\tau)=\chi(\tau/L)\widehat {1_T}(\tau)$$
for some $L >0$.
\begin{lem} \label{cutoff lema}
For any $L, T>0$, it holds
\begin{align}
 \|1^{high}_{T,L} \|_{L^{\frac{3}{2}}(\Real)}\lesssim T^{\frac{1}{3}}L^{-\frac{1}{3}}, \label{estimate571}
\end{align}
\begin{align}
  \|1^{high}_{T,L}\|_{L^{\infty}}\lesssim 1, \label{estimate572}
\end{align}
\begin{align}
  \|1^{low}_{T,L}\|_{L^{\infty}}\lesssim 1. \label{estimate573}
\end{align}
\end{lem}
{\bf Proof.} A direct computation gives
 \begin{align}
\|1^{high}_{T,L} \|_{L^{\frac{3}{2}}(\Real)}&=(\int_\Real|\int_\Real (1_T(t)-1_T(t-\frac{s}{L}))\F^{-1}\chi(s)ds|^{\frac{3}{2}}dt)^{\frac{2}{3}} \notag \\
&\leq\int_\Real(\int_\Real |1_T(t)-1_T(t-\frac{s}{L})|^{\frac{3}{2}}dt)^{\frac{2}{3}}|\F^{-1}\chi(s)|ds \notag \\
&\leq\int_\Real(\int_{[0,T] \backslash[\frac{s}{L},T+\frac{s}{L}]\cup[\frac{s}{L},T+\frac{s}{L}]\backslash[0,T]} dt)^{\frac{2}{3}}|\F^{-1}\chi(s)|ds \notag \\
&\leq\int_\Real(T\wedge\frac{|s|}{L})^{\frac{2}{3}}|\F^{-1}\chi(s)|ds \notag \\
&\lesssim T^{\frac{2}{3}}\wedge L^{-\frac{2}{3}} \notag \\
&\lesssim T^{\frac{1}{3}}L^{-\frac{1}{3}}. \nonumber
\end{align}
The proof of \eqref{estimate572} and \eqref{estimate573} are obvious.

We start with estimates on dyadic pieces. For a dyadic number $N$ and a smooth function $u\in C^\infty_0(\Real^3)$, we write $u_N=P_Nu$.
\begin{lem} \label{L2U2 lema}
Let $N_1, N_2$ be dyadic numbers and $N_2\lesssim N_1$. We have
 \\$(i)$ if $\Omega_1:=\{\lambda=\lambda_1+\lambda_2,|\eta|\leq|\xi|,|\xi|\lesssim |\xi_1-\xi_2|^{\frac{1}{2}}|\xi_1+\xi_2|^{\frac{1}{2}}\}$, then
\begin{align}
 \|\int_{\Omega_1}(\xi+\eta)\widehat {u_{N_1}}(\lambda_1)\widehat {v_{N_2}}(\lambda_2)d\lambda_1\|_{L^2_{\lambda}(\Real^3)}\lesssim N_2^{\frac{1}{2}}\|u_{N_1}\|_{U^2_S}\|v_{N_2}\|_{U^2_S}, \label{estimate574}
\end{align}
 $(ii)$  if $\Omega_2:=\{\lambda=\lambda_1+\lambda_2,|\eta|\leq|\xi|,|\xi|\sim |\xi_1|\sim|\xi_2|\lesssim|\eta_1|\sim|\eta_2|\}$, then
\begin{align}
 \|\int_{\Omega_2}(\xi+\eta)\widehat {u_{N_1}}(\lambda_1)\widehat {v_{N_2}}(\lambda_2)d\lambda_1\|_{L^2_{\lambda}(\Real^3)}\lesssim N_2^{\frac{1}{2}}\|u_{N_1}\|_{U^2_S}\|v_{N_2}\|_{U^2_S}. \label{estimate575}
\end{align}
\end{lem}
{\bf Proof.}  From Proposition \ref{UV prop2}, it suffices to prove the estimates for free solutions $u(t)=e^{tS}u_0$ and $v(t)=e^{tS}v_0$
with $\|u_{N_1}\|_{U^2_S}\|v_{N_2}\|_{U^2_S}$ replaced by $\|P_{N_1}u_0\|_{L^2_{x,y}}\|P_{N_2}v_0\|_{L^2_{x,y}}$.

We follow the argument in the proof of Proposition 1 of \cite{GH14}.

One may assume that $\widehat {u_0}, \widehat {v_0} \geq0$.
\begin{align}
& \int_{\Omega_1}(\xi+\eta)\widehat {e^{tS}u_{0,N_1}}(\lambda_1)\widehat {e^{tS}v_{0,N_2}}(\lambda_2)d\lambda_1\notag \\
=&\int_{\Omega^{'}_1}(\xi+\eta)\delta(\tau-\xi^3_1-\eta^3_1-\xi^3_2-\eta^3_2)\widehat {u_{0,N_1}}(\xi_1,\eta_1)\widehat {v_{0,N_2}}(\xi_2,\eta_2)d\xi_1d\eta_1\notag \\
  =&\frac{1}{3}\int_{\Omega^{''}_1}(\xi+\eta)|\xi(\xi^{*}_1-\xi^{*}_2)|^{-1}(\widehat {u_{0,N_1}}(\xi^{*}_1,\eta_1)\widehat {v_{0,N_2}}(\xi^{*}_2,\eta_2)+\widehat {u_{0,N_1}}
     (\xi^{*}_2,\eta_1)\widehat {v_{0,N_2}}(\xi^{*}_1,\eta_2))d\eta_1, \nonumber
\end{align}
where $\xi^{*}_1$ and $\xi^{*}_2$ are two roots of $g(\xi_1)=\tau-\xi^3_1-\eta^3_1-(\xi-\xi_1)^3-\eta^3_2=0$,
$$\Omega^{'}_1:=\{\zeta=\zeta_1+\zeta_2,|\eta|\leq|\xi|,|\xi|\lesssim |\xi_1-\xi_2|^{\frac{1}{2}}|\xi_1+\xi_2|^{\frac{1}{2}}\}\subset\Real^2,$$
$$\Omega^{''}_1:=\{\eta=\eta_1+\eta_2,|\eta|\leq|\xi|,|\xi|\lesssim |\xi^{*}_1-\xi^{*}_2|^{\frac{1}{2}}|\xi^{*}_1+\xi^{*}_2|^{\frac{1}{2}}\}\subset \Real.$$
Therefore, \eqref{estimate54} gives
\begin{align}
& \|\int_{\Omega_1}(\xi+\eta)\widehat {e^{tS}u_{0,N_1}}(\lambda_1)\widehat {e^{tS}v_{0,N_2}}(\lambda_2)d\lambda_1\|_{L^2_{\lambda}(\Real^3)}\notag \\
  \lesssim &\|\int_{\Omega^{'}_1}\frac{|\xi+\eta|}{|\xi(\xi_1-\xi_2)|^{\frac{1}{2}}}|\xi(\xi_1-\xi_2)|^{\frac{1}{2}}\delta(\tau-\xi^3_1-\eta^3_1-\xi^3_2-\eta^3_2)\widehat {u_{0,N_1}}(\xi_1,\eta_1)\widehat {v_{0,N_2}}(\xi_2,\eta_2)d\xi_1d\eta_1\|_{L^2_{\lambda}(\Real^3)}\notag \\
  \lesssim &\|\int_{\Omega^{'}_1}|\xi(\xi_1-\xi_2)|^{\frac{1}{2}}\delta(\tau-\xi^3_1-\eta^3_1-\xi^3_2-\eta^3_2)\widehat {u_{0,N_1}}(\xi_1,\eta_1)\widehat {v_{0,N_2}}(\xi_2,\eta_2)d\xi_1d\eta_1\|_{L^2_{\lambda}(\Real^3)}\notag \\
  \lesssim &\|\F I^{\frac{1}{2}}_xI^{\frac{1}{2}}_{x,-}(e^{tS}u_{0,N_1},e^{tS}v_{0,N_2})\|_{L^2_{\lambda}(\Real^3)}\notag \\
  \lesssim & N_2^{\frac{1}{2}}\|u_{0,N_1}\|_{L^2_{x,y}}\|v_{0,N_2}\|_{L^2_{x,y}}.  \nonumber
\end{align}

In addition, we have $|\xi|\lesssim N_2^{\frac{1}{2}}|\xi_1|^{\frac{1}{8}}|\eta_1|^{\frac{1}{8}}|\xi_2|^{\frac{1}{8}}|\eta_2|^{\frac{1}{8}}$ on $\Omega_2$.

Denote $\Omega^{'}_2=\{\zeta=\zeta_1+\zeta_2,|\eta|\leq|\xi|,|\xi|\sim |\xi_1|\sim|\xi_2|\lesssim|\eta_1|\sim|\eta_2|\}\subset\Real^2$. From H\"{o}rder inequalities and
\eqref{estimate53}, one gets
\begin{align}
& \|\int_{\Omega_2}(\xi+\eta)\widehat {e^{tS}u_{0,N_1}}(\lambda_1)\widehat {e^{tS}v_{0,N_2}}(\lambda_2)d\lambda_1\|_{L^2_{\lambda}(\Real^3)}\notag \\
  \lesssim &\|\int_{\Omega^{'}_2}|\xi|\delta(\tau-\xi^3_1-\eta^3_1-\xi^3_2-\eta^3_2)\widehat {u_{0,N_1}}(\xi_1,\eta_1)\widehat {v_{0,N_2}}(\xi_2,\eta_2)
             d\xi_1d\eta_1\|_{L^2_{\lambda}(\Real^3)}\notag \\
  \lesssim & N_2^{\frac{1}{2}}\|\int_{\Omega^{'}_2}|\xi_1|^{\frac{1}{8}}|\eta_1|^{\frac{1}{8}}|\xi_2|^{\frac{1}{8}}|\eta_2|^{\frac{1}{8}}\delta(\tau-\xi^3_1-\eta^3_1-\xi^3_2-\eta^3_2)
             \widehat {u_{0,N_1}}(\xi_1,\eta_1)\widehat {v_{0,N_2}}(\xi_2,\eta_2)d\xi_1d\eta_1\|_{L^2_{\lambda}(\Real^3)}\notag \\
  \lesssim & N_2^{\frac{1}{2}}\|\int_{\Omega_2}\F({I^{\frac{1}{8}}_xI^{\frac{1}{8}}_ye^{tS}u_{0,N_1}})(\lambda_1)\F({I^{\frac{1}{8}}_xI^{\frac{1}{8}}_ye^{tS}v_{0,N_2}})(\lambda_2)
             d\lambda_1\|_{L^2_{\lambda}(\Real^3)}\notag \\
  \lesssim & N_2^{\frac{1}{2}}\|I^{\frac{1}{8}}_xI^{\frac{1}{8}}_ye^{tS}u_{0,N_1}\|_{L^4(\Real^3)}\|I^{\frac{1}{8}}_xI^{\frac{1}{8}}_ye^{tS}v_{0,N_2}\|_{L^4(\Real^3)}\notag \\
  \lesssim & N_2^{\frac{1}{2}}\|u_{0,N_1}\|_{L^2_{x,y}}\|v_{0,N_2}\|_{L^2_{x,y}}.  \nonumber
\end{align}

We complete the proof of this lemma.
\begin{prop} \label {Bilinear0}  Let $N_1, N_2, N_3$ be dyadic numbers. There exists $C>0$ such that for all $0<T\leq1$ it holds
\\if $N_1\sim N_3\gg N_2$,
\begin{align}
&|\iiint_{\Real^3}\chi(\frac{t}{T})u_{N_1}v_{N_2}(\partial_x+\partial_y)w_{N_3}dxdydt|\notag \\
\leq& CT^{\frac{1}{2}}N_2^{\frac{1}{2}}\|u_{N_1}\|_{U^2_S}\|v_{N_2}\|_{U^2_S}\|w_{N_3}\|_{V^2_S} \label{estimate55},\end{align}
if $N_1\sim  N_2\gtrsim N_3$,
\begin{align}
&|\iiint_{\Real^3}\chi(\frac{t}{T})u_{N_1}v_{N_2}(\partial_x+\partial_y)w_{N_3}dxdydt|\notag \\
\leq& CT^{\frac{1}{6}}N_1^{\frac{1}{2}}\|u_{N_1}\|_{U^2_S}\|v_{N_2}\|_{U^2_S}\|w_{N_3}\|_{V^2_S}.\label{estimate56}
\end{align}
\end{prop}
{\bf Proof.} It follows from Parseval formula that
\begin{align}
&\iiint_{\Real^3}\chi(\frac{t}{T})u_{N_1}v_{N_2}(\partial_x+\partial_y)w_{N_3}dxdydt\notag \\
=&\int_{\sum^3_{j=1}\lambda_j=0}(\xi_3+\eta_3)\widehat {u_{N_1}}(\lambda_1)\widehat {v_{N_2}}(\lambda_2)\widehat {\chi(\frac{t}{T})w_{N_3}}(\lambda_3).\nonumber
\end{align}
One can assume by symmetry that $|\eta_3|\leq|\xi_3|$ and denote
$$R=\int_{\sum^3_{j=1}\lambda_j=0}(\xi_3+\eta_3)\widehat {u_{N_1}}(\lambda_1)\widehat {v_{N_2}}(\lambda_2)\widehat {\chi(\frac{t}{T})w_{N_3}}(\lambda_3).$$

If $N_1\sim N_3\gg N_2$, then $|\xi_3|\lesssim |\xi_1-\xi_2|$. Otherwise, $|\xi_3|\gg |\xi_1-\xi_2|$. Notice that $\xi_1+\xi_2+\xi_3=0$, we have
$$|\xi_1|\sim|\xi_2|\sim|\xi_3|\lesssim N_2\ll N_1,$$
so
$$|\eta_1|\sim N_1 \gg N_2 \geqslant |\eta_2|.$$
However $\eta_1+\eta_2+\eta_3=0$ implies
$$|\eta_3|\sim |\eta_1|\sim N_1 \gg N_2\gtrsim |\xi_3|$$
which is contradicted against our assumption.

Therefore, we have $|\xi_3|\lesssim |\xi_1-\xi_2|^{\frac{1}{2}}|\xi_1+\xi_2|^{\frac{1}{2}}$ under such circumstance.
Applying H\"{o}lder inequalities, estimates \eqref {estimate574} and \eqref {estimate51}, we can obtain
\begin{align}
R&\lesssim \|\int_{\Omega_1}(\xi+\eta)\widehat {u_{N_1}}(\lambda_1)\widehat {v_{N_2}}(\lambda_2)d\lambda_1\|_{L^2(\Real^3)}\|\chi(\frac{t}{T})w_{N_3}\|_{L^2(\Real^3)}\notag \\
&\lesssim T^{\frac{1}{2}}N_2^{\frac{1}{2}}\|u_{N_1}\|_{U^2_S}\|v_{N_2}\|_{U^2_S}\|w_{N_3}\|_{V^2_S}. \nonumber
\end{align}

We split the domain of the integration into five regions so as to prove \eqref{estimate56}. $R=R_1+R_2+R_3+R_4+R_5.$ Because of the same status of $u$ and $v$, one can assume $|\eta_1|\geq |\eta_2|$.
\\{\bf Region 1.} $|\xi_3|\lesssim |\xi_1-\xi_2|$

Using the same trick as above, then
$$R_1\lesssim T^{\frac{1}{2}}N_1^{\frac{1}{2}}\|u_{N_1}\|_{U^2_S}\|v_{N_2}\|_{U^2_S}\|w_{N_3}\|_{V^2_S}.$$
\\{\bf Region 2.} $|\xi_3|\gg |\xi_1-\xi_2|$ and $|\eta_1|\gg |\xi_3|$

Under this condition, we have
$$|\eta_1|\sim |\eta_2|\gg |\xi_3|\sim|\xi_1|\sim|\xi_2|.$$
Applying H\"{o}lder inequalities, \eqref{estimate575} and \eqref{estimate51}, we obtain
\begin{align}
R_2&\lesssim  \|\int_{\Omega_2}(\xi+\eta)\widehat {u_{N_1}}(\lambda_1)\widehat {v_{N_2}}(\lambda_2)d\lambda_1\|_{L^2(\Real^3)} \|\chi(\frac{t}{T})w_{N_3}\|_{L^2(\Real^3)}\notag \\
&\lesssim T^{\frac{1}{2}}N_1^{\frac{1}{2}}\|u_{N_1}\|_{U^2_S}\|v_{N_2}\|_{U^2_S}\|w_{N_3}\|_{V^2_S}.\nonumber
\end{align}
{\bf Region 3.} $|\xi_3|\gg |\xi_1-\xi_2|$ and $|\eta_1|\sim |\xi_3|\sim |\eta_2|$

We can take the same technique as in Region 2 to get the bound of $R_3$.
\\{\bf Region 4.} $|\xi_3|\gg |\xi_1-\xi_2|$ and $|\eta_1|\sim |\xi_3|\gg |\eta_2|$

Note that $\eta_1+\eta_2+\eta_3=0$, so
$$|\eta_1|\sim |\eta_3|\sim|\xi_3|,$$
hence
$$|\xi_1\xi_2\xi_3|\gg|\eta_1\eta_2\eta_3|.$$

We decompose $Id=Q^S_{<M}+Q^S_{\geq M}$, where $M$ will be chosen later, and we divide the integral $R_4$ into eight pieces of the form
$$\int_{\sum^3_{j=1}\lambda_j=0}(\xi_3+\eta_3)\widehat {1_TQ^S_1u_{N_1}}(\lambda_1)\widehat {1_TQ^S_2v_{N_2}}(\lambda_2)\widehat {1_T
Q^S_3w_{N_3}}(\lambda_3)$$
with $Q^S_j\in \left\{Q^S_{<M},Q^S_{\geq M}\right\}$, $j=1,2,3$.
\\{\bf Case A.} $Q^S_j=Q^S_{<M}$ for $j=1,2,3$

We go a step further to decompose time cut-off as low- and high-frequency parts.
\\{\bf Case A(1).} All of these three are low-frequence

That's to say we need to estimate
\begin{align}
&\int_{\sum^3_{j=1}\lambda_j=0}(\xi_3+\eta_3)\widehat {1^{low}_{T,L}Q^S_{<M}u_{N_1}}(\lambda_1)\widehat {1^{low}_{T,L}Q^S_{<M}v_{N_2}}(\lambda_2)\widehat {1^{low}_{T,L}Q^S_{<M}w_{N_3}}(\lambda_3) \notag \\
=&\int_*(\xi_3+\eta_3)\widehat {1^{low}_{T,L}}(\tau_4)\widehat {Q^S_{<M}u_{N_1}}(\lambda_1)\widehat {1^{low}_{T,L}}(\tau_5)\widehat {Q^S_{<M}v_{N_2}}(\lambda_2)\widehat {1^{low}_{T,L}}(\tau_6)\widehat {Q^S_{<M}w_{N_3}}(\lambda_3).\nonumber
\end{align}
We have that $|\tau_j|<L$  for $j=4,5,6$ and $|\mu_j|<M$  for $j=1,2,3$ due to the cut off operators $1^{low}_{T,L}$ and $Q^S_{<M}$, where $\mu_j=\tau_j-\xi_j^3-\eta_j^3$.

The hyperplane $\left\{\sum^6_{j=1}\tau_j=0,\sum^3_{j=1}\xi_j=0,\sum^3_{j=1}\eta_j=0\right\}$ implies
$$|\mu_1+\mu_2+\mu_3|=|\tau_4+\tau_5+\tau_6+3(\xi_1\xi_2\xi_3+\eta_1\eta_2\eta_3)|.$$
We choose $L=\frac{1}{1000}N_1^2N_3\ll|\xi_1\xi_2\xi_3|$. Thus,
$$M>max(|\mu_1|,|\mu_2|,|\mu_3|)\geq\frac{1}{2}|\xi_1\xi_2\xi_3|>\frac{1}{100}N_1^2N_3$$
holds within the domain of integration.

Therefore, if we set $M=\frac{1}{100}N_1^2N_3$, it follows that
$$\int_{\sum^3_{j=1}\lambda_j=0}(\xi_3+\eta_3)\widehat {1^{low}_{T,L}Q^S_{<M}u_{N_1}}(\lambda_1)\widehat {1^{low}_{T,L}Q^S_{<M}v_{N_2}}(\lambda_2)\widehat {1^{low}_{T,L}Q^S_{<M}w_{N_3}}(\lambda_3) =0.$$
\\{\bf Case A(2).} At least one of these three is high-frequence

For example, we give the estimate when the first one is high-frequence
$$\int_{\sum^3_{j=1}\lambda_j=0}(\xi_3+\eta_3)\widehat {1^{high}_{T,L}Q^S_{<M}u_{N_1}}(\lambda_1)\widehat {1_vQ^S_{<M}v_{N_2}}(\lambda_2)\widehat {1_wQ^S_{<M}w_{N_3}}(\lambda_3),$$
where $1_v,1_w\in \left\{1^{high}_{T,L},1^{low}_{T,L}\right\}$.

H\"{o}lder inequalities, \eqref{estimate52}, Lemma \ref{cutoff lema} and \eqref {UV estimate 2} provide
\begin{align}
&|\int_{\sum^3_{j=1}\lambda_j=0}(\xi_3+\eta_3)\widehat {1^{high}_{T,L}Q^S_{<M}u_{N_1}}(\lambda_1)\widehat {1_vQ^S_{<M}v_{N_2}}(\lambda_2)\widehat {1_wQ^S_{<M}w_{N_3}}(\lambda_3)| \notag \\
\lesssim& N_3\|1^{high}_{T,L}Q^S_{<M}u_{N_1}\|_{L^{\frac{9}{7}}_tL^3_{xy}}\|1_vQ^S_{<M}v_{N_2}\|_{L^9_tL^3_{xy}}\|1_wQ^S_{<M}w_{N_3}\|_{L^9_tL^3_{xy}}    \notag \\
\lesssim& N_3\|1^{high}_{T,L}\|_{L^{\frac{3}{2}}_t}\|Q^S_{<M}u_{N_1}\|_{L^9_tL^3_{xy}}\|1_v\|_{L^{\infty}_t}\|Q^S_{<M}v_{N_2}\|_{L^9_tL^3_{xy}}\|1_w\|_{L^{\infty}_t}
    \|Q^S_{<M}w_{N_3}\|_{L^9_tL^3_{xy}}   \notag \\
\lesssim& N_3T^{\frac{1}{3}}N^{-\frac{1}{3}}_3N^{-\frac{2}{3}}_1\|Q^S_{<M}u_{N_1}\|_{U^9_S}\|Q^S_{<M}v_{N_2}\|_{U^9_S}\|Q^S_{<M}w_{N_3}\|_{U^9_S} \notag \\
\lesssim& T^{\frac{1}{3}}N^{\frac{2}{3}}_1N^{-\frac{2}{3}}_1\|u_{N_1}\|_{U^9_S}\|v_{N_2}\|_{U^9_S}\|w_{N_3}\|_{U^9_S} \notag \\
\lesssim& T^{\frac{1}{3}}N^{\frac{1}{2}}_1\|u_{N_1}\|_{U^2_S}\|v_{N_2}\|_{U^2_S}\|w_{N_3}\|_{V^2_S}.\nonumber
\end{align}
{\bf Case B.} $Q^S_j=Q^S_{\geq M}$ for some $j=1,2,3$

We take $Q^S_3=Q^S_{\geq M}$ for instance to bound $R_5$ by the right hand of \eqref{estimate56}. The cases $j=1,2$ can be dealt with in the same way.

Using H\"{o}lder inequalities, \eqref {UV estimate 1}, \eqref{estimate52} and \eqref {UV estimate 2}, we have
\begin{align}
&|\int_{\sum^3_{j=1}\lambda_j=0}(\xi_3+\eta_3)\widehat {\chi(\frac{t}{T})Q^S_1u_{N_1}}(\lambda_1)\widehat {\chi(\frac{t}{T})Q^S_2v_{N_2}}(\lambda_2)
\widehat {\chi(\frac{t}{T})Q^S_{\geq M}w_{N_3}}(\lambda_3)|\notag \\
\lesssim& N_3\|\chi(\frac{t}{T})Q^S_1u_{N_1}\|_{L^4(\Real^3)}\|\chi(\frac{t}{T})Q^S_2v_{N_2}\|_{L^4(\Real^3)}\|Q^S_{\geq M}w_{N_3}\|_{L^2(\Real^3)} \notag \\
\lesssim& N_3M^{-\frac{1}{2}}\|\chi(\frac{t}{T})\|^2_{L^{12}(\Real)}\|Q^S_1u_{N_1}\|_{L^6_tL^4_{x,y}}\|Q^S_2v_{N_2}\|_{L^6_tL^4_{x,y}}
\|w_{N_3}\|_{V^2_S}\notag \\
\lesssim& T^{\frac{1}{6}}N_1^{-1}N_3^{\frac{1}{2}}\|Q^S_1u_{N_1}\|_{U^6_S}\|Q^S_2v_{N_2}\|_{U^6_S}\|w_{N_3}\|_{V^2_S}\notag \\
\lesssim& T^{\frac{1}{6}}N_1^{-1}N_1^{\frac{1}{2}}\|u_{N_1}\|_{U^6_S}\|v_{N_2}\|_{U^6_S}\|w_{N_3}\|_{V^2_S}\notag \\
\lesssim& T^{\frac{1}{6}}N_1^{\frac{1}{2}}\|u_{N_1}\|_{U^2_S}\|v_{N_2}\|_{U^2_S}\|w_{N_3}\|_{V^2_S}.\nonumber
\end{align}
{\bf Region 5.} $|\xi_3|\gg |\xi_1-\xi_2|$ and $|\xi_3|\gg |\eta_1|\geq |\eta_2|$

Because $|\xi_1\xi_2\xi_3|\gg|\eta_1\eta_2\eta_3|$ also holds true in this case, one can estimate $R_5$ just like $R_4$.

Hence, the proof of this proposition is complete. \vspace{3mm}

We denote by $Y^s$ the space of all function $u\in \mathscr{S}'(\Real^3)$ such that
$$\|u\|_{Y^s}:=\sum_N N^s\|P_Nu\|_{U^2_S}<\infty.$$
The work space we choose is $Y^{\frac{1}{2}}$.
\begin{prop} \label {Bilinear1} Let $0<T\leq1$. We have
\begin{align}
&\|\int^t_0e^{(t-t')S}\chi(\frac{t}{T})(\partial_x+\partial_y)(uv)(t')dt'\|_{Y^{\frac{1}{2}}}\lesssim T^{\frac{1}{6}}\|u\|_{Y^{\frac{1}{2}}}\|v\|_{Y^{\frac{1}{2}}} \label{estimate57}.
 \end{align}
\end{prop}
{\bf Proof.} From symmetry and definition of $Y^{\frac{1}{2}}$ we need to consider the two terms
\begin{align}
J_1:=\sum_{N}\sum_{N_1\gg N_2}N^{\frac{1}{2}}\|\int^t_0e^{(t-t')S}\chi(\frac{t}{T})(\partial_x+\partial_y)P_N(u_{N_1}v_{N_2})(t')dt'\|_{U^2_S}\nonumber
\end{align}
and
\begin{align}
J_2:=\sum_{N}\sum_{N_1\sim N_2}N^{\frac{1}{2}}\|\int^t_0e^{(t-t')S}\chi(\frac{t}{T})(\partial_x+\partial_y)P_N(u_{N_1}v_{N_2})(t')dt'\|_{U^2_S}.\nonumber
\end{align}
By Proposition \ref{UV prop1} \rm(iii) and \eqref{estimate55}, one has
\begin{align}
J_1&\lesssim \sum_{N}\sum_{N_1\gg N_2}N^{\frac{1}{2}}\mathop{sup}\limits_{\|w\|_{V^2_S}\leq1}|\iiint_{\Real^3}\chi(\frac{t}{T})u_{N_1}v_{N_2}
(\partial_x+\partial_y)w_{N}dxdydt|   \notag \\
&\lesssim T^{\frac{1}{2}}\sum_{N_1}\sum_{N_1\gg N_2}N_1^{\frac{1}{2}}N_2^{\frac{1}{2}}\|u_{N_1}\|_{U^2_S}\|v_{N_2}\|_{U^2_S}    \notag \\
&\lesssim T^{\frac{1}{2}}\|u\|_{Y^{\frac{1}{2}}}\|v\|_{Y^{\frac{1}{2}}}.\nonumber
 \end{align}
Whereas the second one can be controlled with the help of \eqref{estimate56},
\begin{align}
J_2&\lesssim \sum_{N}\sum_{N_1\sim N_2}N^{\frac{1}{2}}\mathop{sup}\limits_{\|w\|_{V^2_S}\leq1}|\iiint_{\Real^3}\chi(\frac{t}{T})u_{N_1}v_{N_2}
(\partial_x+\partial_y)w_{N}dxdydt|   \notag \\
   &\lesssim T^{\frac{1}{6}}\sum_{N\lesssim N_1}\sum_{N_1}N^{\frac{1}{2}}N_1^{\frac{1}{2}}\|u_{N_1}\|_{U^2_S}\|v_{N_1}\|_{U^2_S}    \notag \\
  &\lesssim T^{\frac{1}{6}}\sum_{N_1}N_1\|u_{N_1}\|_{U^2_S}\|v_{N_1}\|_{U^2_S}    \notag \\
&\lesssim T^{\frac{1}{6}}\|u\|_{Y^{\frac{1}{2}}}\|v\|_{Y^{\frac{1}{2}}}.\nonumber
 \end{align}

\section{Local well-posedness}
For the sake of completeness, we sketch the proof of Theorem 1.2 based on \eqref{estimate57}. We rewrite \eqref{ZK1} as an integral equation
$$u=\mathscr{T}u$$
where
$$\mathscr{T}u:=\chi(\frac{t}{T})e^{tS}u_0-\int^t_0e^{(t-t')S}\chi(\frac{t}{T})(\partial_x+\partial_y)(u^2)(t')dt'.$$
Like inequality (6.6) in \cite{MP15}, we can get
$$\|\chi(\frac{t}{T})e^{tS}u_0\|_{Y^{\frac{1}{2}}}\lesssim \|u_0\|_{B^{\frac{1}{2}}_{2,1}}.$$
Proposition \ref {Bilinear1} gives us that
$$\|\int^t_0e^{(t-t')S}\chi(\frac{t}{T})(\partial_x+\partial_y)(u^2)(t')dt'\|_{Y^{\frac{1}{2}}}
\lesssim T^{\frac{1}{6}}\|u\|^2_{Y^{\frac{1}{2}}}.$$
Hence, we have
$$\|\mathscr{T}u\|_{Y^{\frac{1}{2}}}<C_0(\|u_0\|_{B^{\frac{1}{2}}_{2,1}}+T^{\frac{1}{6}}\|u\|^2_{Y^{\frac{1}{2}}}),$$
and similarly
$$\|\mathscr{T}u-\mathscr{T}v\|_{Y^{\frac{1}{2}}}<C_0T^{\frac{1}{6}}\|u-v\|_{Y^{\frac{1}{2}}}\|u+v\|_{Y^{\frac{1}{2}}}.$$

Defining
$$B_r:=\{u\in Y^{\frac{1}{2}}_T\ | \ \ \|u\|_{Y^{\frac{1}{2}}_T}<r\}$$
with $r=4C_0\|u_0\|_{B^{\frac{1}{2}}_{2,1}}$ and $T=min\left\{1, \frac{1}{(4C_0r)^6}\right\}$, it is easy to verify that
$$\mathscr{T}: B_r\rightarrow B_r$$
is a strict contraction. Therefore there exists a unique fixed point in $B_r$, which solves \eqref{ZK1} on the interval $[0,T]$.\vspace{4mm}
\\{\bf Acknowledgment}\vspace{2mm}
\\The author would like to express his deep gratitude to Professor Yoshio Tsutsumi for illuminating discussion as well as Professor Nobu Kishimoto for sharing him with the time cut-off decomposition technique which is introduced in the paper of Molinet and Vento. He  is glad to acknowledge his indebtedness for the suggestion he received from Professor Luc Molinet and Professor Baoxiang Wang. He also thanks the anonymous referees for all kinds of significant advice.

\nocite{*}
\bibliography{document}

\footnotesize

\begin{thebibliography}{99}

\bibitem{BKS03} Matania Ben-Artzi, Herbert Koch, and Jean-Claude Saut. Dispersion estimates for third order equations in two dimensions. Comm. Partial Differential Equations, \textbf{28}(11-12) (2003), 1943--1974.

\bibitem{Bourgain 2a} J. Bourgain. Exponential sums and nonlinear Schrödinger equations. Geom. Funct. Anal., \textbf{3}(2) (1993), 157--178.

\bibitem{CKSTT02} J. Colliander, M. Keel, G. Staffilani, H. Takaoka, and T. Tao. Almost conservation laws and global rough solutions to a Nonlinear Schr\"{o}dinger equation. Math. Res. Lett., \textbf{9}(5-6) (2002), 659--682.

\bibitem{CKSTT04} J. Colliander, M. Keel, G. Staffilani, H. Takaoka, and T. Tao. Multilinear estimates for periodic KdV equations, and applications. J. Funct. Anal., \textbf{211}(1) (2004), 173--218.

\bibitem{F95} Andrei V. Faminskii, The Cauchy problem for the Zakharov-Kuznetsov equation. Differ. Equations, \textbf{31}(6) (1995), 1002--1012.

\bibitem{GTV97} Jean Ginibre, Yoshio Tsutsumi, and Giorgio Velo. On the Cauchy problem for the Zakharov system. J. Funct. Anal., \textbf{151}(2) (1997), 384--436.

\bibitem{GH14} Axel Gr\"{u}nrock and Sebastian Herr. The Fourier restriction norm method for the Zakharov-Kuznetsov equation. Discrete Contin. Dyn. Syst., \textbf{34}(5) (2014), 2061--2068.

\bibitem{HaHeKo09} M. Hadac, S. Herr and H. Koch. Well-posedness and scattering for the KP-II equation in a critical space. Ann. Inst. H. Poincar\'{e} Anal. Non
Lin\'{e}aire, 26 (2009),  917-941.

\bibitem{KPV93} Carlos E. Kenig, Gustavo Ponce, and  Luis Vega. Well-posedness and scattering results for the generalized Korteweg-de Vries equation via the contraction principle. Comm. Pure Appl. Math., \textbf{46}(4) (1993), 527--620.

\bibitem{KoTa05} H. Koch and D. Tataru. Dispersive estimates for principlally normal pseudo-differential operators, Comm. Pure Appl. Math., \textbf{58}(2005), 217--284.

\bibitem{KoTa07} H. Koch and D. Tataru. A priori bounds for the 1D cubic NLS in negative Sobolev spaces, Int. Math. Res. Not., 2007, no. 16, Art.
ID rnm053, 36 pp.

\bibitem{KoTa12} H. Koch and D. Tataru. Energy and local energy bounds for the 1D cubic NLS equation in $H^{\frac{1}{4}}$,  Ann. Inst. H. Poincar\'{e} Anal. Non
Lin\'{e}aire, \textbf{29} (2012),  955--988.

\bibitem{KT03} Herbert Koch and Nikolay Tzvetkov. On the local well-posedness of the Benjamin-Ono equation in $H^s(\Real)$. Int. Math. Res. Not., \textbf{26} (2003), 1449--1464.

\bibitem{ZK74} E.A. Kuznetsov and V.E. Zakharov. On the three dimensional solitons. Sov. Phys. JETP, \textbf{39} (1974), 285--286.

\bibitem{EWLKHS82} Ernst W. Laedke and Karl-Heinz Spatschek. Nonlinear ion-acoustic waves in weak magnetic fields. Phy. Fluids, \textbf{25}(6) (1982), 985--989.

\bibitem{LLS13} David Lannes, Felipe Linares, and Jean-Claude Saut. The Cauchy problem for the Euler-Poisson system and derivation of the Zakharov-Kuznetsov equation. Studies in phase space analysis with applications to PDEs, Springer New York (2013), 181--213.

\bibitem{LP09} Felipe Linares and Ademir Pastor. Well-posedness for the two-dimensional modified Zakharov-Kuznetsov equation. SIAM J. Math. Anal., \textbf{41}(4) (2009), 1323--1339.

\bibitem{FLJCS82} Felipe Linares and Jean-Claude Saut. The Cauchy problem for the 3D Zakharov-Kuznetsov equation. Discrete Contin. Dyn. Syst., \textbf{24}(2) (2009), 547--565.

\bibitem{MP15} Luc Molinet and Didier Pilod. Bilinear Strichartz estimates for the Zakharov-Kuznetsov equation and applications. Ann. Inst. H. Poincar\'{e} Anal. Non Lin\'{e}aire, \textbf{32}(2) (2015), 347--371.

\bibitem{MV15} Luc Molinet and St\'{e}phane Vento. Improvement of the energy method for strongly nonresonant dispersive equations and applications. Anal. PDE, \textbf{8}(6) (2015), 1455--1495.

\bibitem{RV12} Felipe Ribaud and St\'{e}phane Vento. Well-Posedness Results for the Three-Dimensional Zakharov-Kuznetsov Equation. SIAM J. Math. Anal.,  \textbf{44}(4) (2012), 2289--2304.

\bibitem{WaHaHuGu11} B. X. Wang,  Z. H. Huo, C. C. Hao and Z. H. Guo, Harmonic Analysis Method for Nonlinear Evolution Equations. I. World Scientific Publishing Co., Pte. Ltd., Hackensack, NJ (2011).

\end{thebibliography}
\bibliographystyle{amsplain}

\end{document}